\pdfoutput=1
\RequirePackage{ifpdf}
\ifpdf 
\documentclass[pdftex]{sigma}
\else
\documentclass{sigma}
\fi

\usepackage[mathscr]{eucal}

\newcommand{\geh}{\mathfrak{g}}
\newcommand{\la}{\lambda}
\newcommand{\ot}{\otimes}
\newcommand{\pair}[1]{\langle{#1}\rangle}

\newcommand{\Z}{{\mathbb Z}}

\newcommand{\C}{{\mathbb C}}
\newcommand{\Q}{{\mathbb Q}}

\def\eshort{e_{1}}
\def\elong{e_{2}}
\def\qmq#1{\langle #1 \rangle}
\def\undividedB#1{{\tilde E}^{#1}_{\bf 2}}

\def\bk{{\bf k}}
\def\bK{{\bf K}}

\def\am{{{\bf a}^-}}

\begin{document}

\allowdisplaybreaks

\renewcommand{\thefootnote}{$\star$}

\renewcommand{\PaperNumber}{049}

\FirstPageHeading

\ShortArticleName{PBW Bases and Quantized Algebra of Functions}

\ArticleName{A Common Structure in PBW Bases \\
of the Nilpotent Subalgebra of $\boldsymbol{U_q(\mathfrak{g})}$\\
and Quantized Algebra of Functions\footnote{This paper is a~contribution to the Special Issue in honor of Anatol Kirillov and Tetsuji Miwa.
The full collection is available at
\href{http://www.emis.de/journals/SIGMA/InfiniteAnalysis2013.html}{http://www.emis.de/journals/SIGMA/InfiniteAnalysis2013.html}}}

\Author{Atsuo KUNIBA~$^\dag$, Masato OKADO~$^\ddag$ and Yasuhiko YAMADA~$^\S$}

\AuthorNameForHeading{A.~Kuniba, M.~Okado and Y.~Yamada}

\Address{$^\dag$~Institute of Physics, Graduate School of Arts and Sciences,
University of Tokyo,\\
\hphantom{$^\dag$}~Komaba, Tokyo 153-8902, Japan}
\EmailD{\href{mailto:atsuo@gokutan.c.u-tokyo.ac.jp}{atsuo@gokutan.c.u-tokyo.ac.jp}}

\Address{$^\ddag$~Department of Mathematical Science,
Graduate School of Engineering Science,\\
\hphantom{$^\ddag$}~Osaka University, Toyonaka, Osaka 560-8531, Japan}
\EmailD{\href{mailto:okado@sigmath.es.osaka-u.ac.jp}{okado@sigmath.es.osaka-u.ac.jp}}

\Address{$^\S$~Department of Mathematics, Faculty of Science,
Kobe University, Hyogo 657-8501, Japan}
\EmailD{\href{mailto:yamaday@math.kobe-u.ac.jp}{yamaday@math.kobe-u.ac.jp}}

\ArticleDates{Received March 19, 2013, in f\/inal form July 10, 2013; Published online July 19, 2013}

\Abstract{For a f\/inite-dimensional simple Lie algebra $\mathfrak{g}$,
let $U^+_q(\mathfrak{g})$ be the positive part of
the quantized universal enveloping algebra, and
$A_q(\mathfrak{g})$ be the quantized algebra of functions.
We show  that the transition matrix of the PBW bases of $U^+_q(\mathfrak{g})$
coincides with the intertwiner between the irreducible
$A_q(\mathfrak{g})$-modules labeled by two dif\/ferent reduced expressions of the longest element of the
Weyl group of $\mathfrak{g}$. This generalizes the earlier result by Sergeev on $A_2$ related to the
tetrahedron equation and endows
a new representation theoretical interpretation with
the recent solution to the 3D ref\/lection equation for $C_2$.
Our proof is based on a realization of $U^+_q(\mathfrak{g})$ in a
quotient ring of $A_q(\mathfrak{g})$.}

\Keywords{quantized enveloping algebra; PBW bases; quantized algebra of functions; tetrahedron equation}

\Classification{17B37; 20G42; 81R50; 17B80}

\begin{flushright}
\begin{minipage}{80mm}
\it Dedicated to Professors Anatol N.~Kirillov and Tetsuji Miwa who taught us the joy of
doing mathematics.
\end{minipage}
\end{flushright}

\renewcommand{\thefootnote}{\arabic{footnote}}
\setcounter{footnote}{0}

\section{Introduction}\label{a:sec1}

Let $\geh$ be a f\/inite-dimensional simple Lie algebra and
$U_q(\geh)$ be the Drinfeld--Jimbo quantized enveloping algebra.
$U_q(\geh)$ has the subalgebra $U^+_q(\geh)$ generated by the
Chevalley generators $e_1,\ldots, e_n$   $(n= \mathrm{rank}\, \geh)$
corresponding to the simple roots.
Denote by $W=\langle s_1,\ldots, s_n\rangle$ the
Weyl group of $\geh$ generated by the simple ref\/lections $s_1,\ldots, s_n$.
It is well known (see for example~\cite{L2}) that
for each reduced expression $w_0 = s_{i_1}\cdots s_{i_l}$
of the longest element $w_0 \in W$,
one can associate the Poincar\'e--Birkhof\/f--Witt (PBW) basis
of $U^+_q(\geh)$ having the form
\begin{gather*}
E_{\bf i}^A=e_{\beta_1}^{(a_1)}e_{\beta_2}^{(a_2)}\cdots e_{\beta_l}^{(a_l)},
\qquad
A=(a_1,\ldots, a_l) \in (\Z_{\ge 0})^l,
\end{gather*}
where $e^{(a_i)}_{\beta_i}$'s are the divided powers of the positive root
vectors determined by the choice ${\bf i}= (i_1,\ldots, i_l)$.
See Section~\ref{a:ss:pbw}.
Let $\big\{E^A_{\bf j}\mid A \in (\Z_{\ge 0})^l\big\}$ with
${\bf j} = (j_1,\ldots, j_l)$ be another
PBW basis associated with a yet dif\/ferent reduced expression
$w_0 = s_{j_1}\cdots s_{j_l}$.
Following Lusztig~\cite{L}, one expands
a basis in terms of another as
\begin{gather*}
E^A_{\bf i} = \sum_{B \in (\Z_{\ge 0})^l}
\gamma^A_B E^B_{\bf j}
\end{gather*}
and obtains the
transition coef\/f\/icient $\gamma^A_B$ uniquely.
We have suppressed its dependence on ${\bf i}$, ${\bf j}$
in the notation.
Many remarkable properties are known for $\gamma^A_B$
including the fact $\gamma^A_B \in \Z[q]$.
See \cite[Proposition~2.3]{L} for example.

In this paper we show that the transition coef\/f\/icients
$\gamma=(\gamma^A_B)$ coincide with
the matrix elements of the intertwiner between the irreducible $A_q(\geh)$-modules labeled by two dif\/ferent reduced expressions of the longest element of the Weyl group of $\geh$.
Here $A_q(\geh)$ denotes the quantized algebra of functions  associated with
$\geh$.
It is a Hopf subalgebra of the dual $U_q(\geh)^\ast$
which has been studied from a variety of aspects.
See \cite{D86,Kas93,NYM,RTF,So1,So2,VS} for example.
Let us brief\/ly recall the most relevant result to the present paper
due to Vaksman and Soibelman \cite{So1,So2,VS}.
To each reduced expression of a (not necessarily longest) element
$w = s_{i_1}\cdots s_{i_r} \in W$,
one can associate an irreducible representation
$\pi_{\bf i}$
labeled by ${\bf i}=(i_1,\ldots, i_r)$ having the form
\begin{gather*}
\pi_{\bf i}= \pi_{i_1}\otimes \cdots \otimes \pi_{i_r}: \
A_q(\geh) \rightarrow
\mathrm{End}\big(\mathcal{F}_{q_{i_1}} \otimes \cdots \otimes
\mathcal{F}_{q_{i_r}}\big),
\end{gather*}
where each component
$\pi_i: A_q(\geh) \rightarrow \mathrm{End}(\mathcal{F}_{q_i})$ is the
{\em fundamental representation} of $A_q(\geh)$ on
the $q$-oscillator Fock space
$\mathcal{F}_{q_i} = \bigoplus_{m\ge 0}\C(q)|m\rangle$.
See Section~\ref{a:ss:gr}.
The two irreducible representations $\pi_{\bf i}$ and
$\pi_{\bf j}$ with ${\bf j}=(j_1,\ldots, j_r)$
are isomorphic if
$s_{i_1}\cdots s_{i_r}=s_{j_1}\cdots s_{j_r} \in W$
are reduced expressions (Theorem \ref{a:th:so}).
Thus one has the intertwiner
$\Phi = \Phi_{{\bf i}, {\bf j}}:
\mathcal{F}_{q_{i_1}}\otimes \cdots \otimes \mathcal{F}_{q_{i_r}}
\rightarrow
\mathcal{F}_{q_{j_1}}\otimes \cdots \otimes \mathcal{F}_{q_{j_r}}$
characterized by
\begin{gather*}
\pi_{\bf j}(g)\circ \Phi = \Phi \circ \pi_{\bf i}(g)\qquad
\forall \, g \in A_q(\mathfrak{g})
\end{gather*}
up to an overall constant.
Writing the basis of the Fock space
$\mathcal{F}_{q_{i_1}} \otimes \cdots \otimes
\mathcal{F}_{q_{i_r}}$ as
$|A\rangle = |a_1\rangle \otimes \cdots \otimes |a_r\rangle$
with $A=(a_1, \ldots, a_r) \in (\Z_{\ge 0})^r$,
we def\/ine the matrix elements of $\Phi=(\Phi^A_B)$ by
$\Phi |B\rangle = \sum_A \Phi^A_B |A\rangle$ and the
normalization $\Phi^{0,\ldots, 0}_{0,\ldots, 0}=1$.
Our main result (Theorem \ref{theorem:gamma=R})
is concerned with the longest element case $r=l$
and is stated for each pair $({\bf i}, {\bf j})$ as $\gamma^A_B = \Phi^A_B$, i.e.,
\begin{gather}\label{a:main}
\gamma= \Phi.
\end{gather}
For a convenience we also introduce the ``checked'' intertwiner
$\Phi^{\vee} = \Phi \circ \sigma$, where
$\sigma  (|a_1\rangle \otimes \cdots \otimes |a_l\rangle) =
 |a_l\rangle \otimes \cdots \otimes |a_1 \rangle$
is the reversal of the components.

Our work is inspired by recent developments in
3-dimensional (3D) integrable systems related to
rank 2 cases.
Recall the Zamolodchikov tetrahedron equation~\cite{Zam80}
and the Isaev--Kulish 3D ref\/lection equation~\cite{IK}:
\begin{gather}
{R}_{356}{R}_{246}{R}_{145}{R}_{123}
={R}_{123}{R}_{145}{R}_{246}{R}_{356}, \label{a:teq}\\
{R}_{456}{R}_{489}
{K}_{3579}{R}_{269}{R}_{258}
{K}_{1678}{K}_{1234}
={K}_{1234}
{K}_{1678}
{R}_{258}{R}_{269}
{K}_{3579}{R}_{489}{R}_{456}.\label{a:3dref}
\end{gather}
They are equalities among the linear operators acting on the
tensor product of~6 and 9 vector spaces, respectively.
The indices specify the components in the tensor product
on which the operators $R$ and $K$ act nontrivially.
They serve as 3D analogue of the Yang--Baxter and ref\/lection equations
postulating certain factorization conditions of straight strings
which undergo the scattering $R$ and the ref\/lection $K$ by a boundary plane.

For $\geh= A_2$,
Kapranov and Voevodsky~\cite{KV} showed that
$R=\Phi^{\vee} \in \mathrm{End}(\mathcal{F}_q^{\otimes 3})$
provides a~solution to the
tetrahedron equation~(\ref{a:teq}).
Moreover it was discovered by Sergeev~\cite{S08}  that
the solution of the tetrahedron equation~$R$ in~\cite{BS} (given also in~\cite{KV} with misprint)
is related with the transition matrix as $\gamma= R \circ \sigma$.
Thus the equality~(\ref{a:main})  for $\geh= A_2$ is a corollary of their results.
Apart from the $A_2$ case, it has been shown more recently~\cite{KO1} that
$K = \Phi^\vee$  for $\geh=C_2$ yields the f\/irst nontrivial
solution to the 3D ref\/lection equation (\ref{a:3dref}).
See also~\cite{KO2} for $\geh=B_2$.
These results motivated us to investigate the general $\geh$ case
and have led to~(\ref{a:main}).
It is our hope that it provides a useful insight into
higher-dimensional integrable systems from
the representation theory of quantum groups.

The layout of the paper is as follows.
In Section~\ref{a:sec2},
we summarize the def\/initions of $U_q(\geh)$ and PBW bases.
In Section~\ref{a:sec3},
we recall  the basic facts on $A_q(\geh)$ following Kashiwara \cite{Kas93}.
A fundamental role is played by the Peter--Weyl type Theorem \ref{th:peter}.
The relation with  the Reshetikhin--Takhtadzhyan--Faddeev realization
by generators and relations~\cite{RTF} is explained and
its concrete forms are quoted for
$A_n$, $C_n$ and~$G_2$~\cite{Sasaki} which will be of use in later sections.
The construction of a certain quotient ring
$A_q(\geh)_{\mathcal S}$ of $A_q(\geh)$ and
the special elements $\sigma_i \in  A_q(\geh)$
(Def\/inition~\ref{def:sigma})
and $\xi_i \in A_q(\geh)_{\mathcal S}$ (\ref{a:xidef}) will play
a key role in
our proof of~(\ref{a:main}).
In Section~\ref{a:sec4},
we brief\/ly review the representation theory  of $A_q(\geh)$
in \cite{So1,So2,VS} and sketch the intertwiners
for the rank~2 cases.
Section~\ref{a:sec5} is devoted
to the proof of the main theorem $\gamma = \Phi$.
It reduces to the rank~2 cases and
is done without recourse to explicit formulae of $\gamma$ or $\Phi$.
Our method is to identify their characterizations
under the correspondence $e_i \mapsto \xi_i$.
Actually, this map extends to an algebra homomorphism
$U^+_q(\geh) \rightarrow A_q(\geh)_{\mathcal S}$
for general $\geh$ as shown by Yakimov~\cite{Yak}.
We give a direct proof of a part of his results in Section~\ref{a:sec6}.

\section[Quantized enveloping algebra $U_q(\geh)$]{Quantized enveloping algebra $\boldsymbol{U_q(\geh)}$}\label{a:sec2}

\subsection{Def\/inition}

In this paper $\geh$ stands for a f\/inite-dimensional simple Lie algebra.
Its weight lattice, simple roots, simple coroots, fundamental weights are
denoted by $P$, $\{\alpha_i\}_{i\in I}$, $\{h_i\}_{i\in I}$, $\{\varpi_i\}_{i\in I}$
where $I$ is the index set of the Dynkin diagram of~$\geh$.
The Cartan matrix $(a_{ij})_{i,j\in I}$ is given by
$a_{ij}=\pair{h_i,\alpha_j}=2(\alpha_i, \alpha_j)/(\alpha_i,\alpha_i)$.

The quantized enveloping algebra $U_q(\geh)$ is an associative algebra over $\Q(q)$ generated by
$\{e_i, f_i, k_i^{\pm1} \mid i\in I\}$ satisfying the relations:
\begin{gather}
 k_i k_j=k_j k_i, \qquad
k_i k_i^{-1}=k_i^{-1}k_i=1, \nonumber
\\
 k_i e_j k_i^{-1}=q_i^{\langle h_i,\alpha_j \rangle} e_j, \qquad
k_i f_j k_i^{-1}=q_i^{-\langle h_i, \alpha_j \rangle} f_j, \qquad
[e_i, f_j]=\delta_{ij}\frac{k_i-k_i^{-1}}{q_i-q_i^{-1}}, \nonumber
\\
 \sum_{r=0}^{1-a_{ij}} (-1)^r e_i^{(r)} e_j e_i^{(1-a_{ij}-r)}
=\sum_{r=0}^{1-a_{ij}} (-1)^r f_i^{(r)} f_j f_i^{(1-a_{ij}-r)}=0,\qquad i\neq j.\label{a:serre}
\end{gather}
Here we use the following notations:
$q_i=q^{(\alpha_i,\alpha_i)/2},[m]_i=(q_i^m-q_i^{-m})/(q_i-q_i^{-1})$,
$[n]_i!=\prod\limits_{m=1}^n[m]_i$, $e_i^{(n)}=e_i^n/[n]_i!,f_i^{(n)}=f_i^n/[n]_i!$.
We normalize the simple roots so that $q_i=q$ when
$\alpha_i$ is a short root.
$U_q(\geh)$ is a Hopf algebra. As its comultiplication we adopt the following one
\begin{gather*}
\Delta(k_i)=k_i\ot k_i,\qquad \Delta(e_i)=e_i\ot1+k_i\ot e_i,\qquad \Delta(f_i)=f_i\ot k_i^{-1}+1\ot f_i.
\end{gather*}

\subsection{PBW basis}\label{a:ss:pbw}

Let $W$ be the Weyl group of $\geh$. It is generated by simple ref\/lections $\{s_i\mid i\in I\}$
obeying the relations: $s_i^2=1$, $(s_is_j)^{m_{ij}}=1$ ($i\ne j$) where $m_{ij}=2,3,4,6$ for
$\pair{h_i,\alpha_j}\pair{h_j,\alpha_i}=0,1,2,3$, respectively. Let $w_0$ be the longest element
of $W$ and f\/ix a reduced expression $w_0=s_{i_1}s_{i_2}\cdots s_{i_l}$. Then every positive root
occurs exactly once in
\[
\beta_1=\alpha_{i_1},\quad \beta_2=s_{i_1}(\alpha_{i_2}), \quad \ldots,\quad
\beta_l=s_{i_1}s_{i_2}\cdots s_{i_{l-1}}(\alpha_{i_l}).
\]
Correspondingly, def\/ine elements $e_{\beta_r}\in U_q(\geh)$ ($r=1,\ldots,l$) by
\begin{gather}\label{a:ebdef}
e_{\beta_r}=T_{i_1}T_{i_2}\cdots T_{i_{r-1}}(e_{i_r}).
\end{gather}
Here $T_i$ is the action of the braid group on $U_q(\geh)$ introduced by Lusztig \cite{L2}. It is
an algebra automorphism and is given on the generators $\{e_j\}$ by
\begin{gather*}
T_i(e_i) = -k_if_i,
\qquad
T_i(e_j)=\sum_{r=0}^{-a_{ij}}(-1)^r q_i^r e_i^{(r)}e_je_i^{(-a_{ij}-r)},
\qquad i\neq j.
\end{gather*}

$U_q(\geh)$ has a subalgebra generated by $\{e_i\mid i\in I\}$, denoted by $U_q^+(\geh)$.
It is known that $e_{\beta_r} \in U^+_q(\geh)$ holds for any $r$.
$U^+_q(\geh)$ has the
so-called Poincar\'e--Birkhof\/f--Witt (PBW) basis. It depends on the reduced expression
$s_{i_1}s_{i_2}\cdots s_{i_l}$ of $w_0$. Set ${\bf i}=(i_1,i_2,\ldots,i_l)$ and def\/ine for
$A=(a_1,a_2,\ldots,a_l)\in(\Z_{\ge0})^l$
\begin{gather}\label{a:Eai}
E_{\bf i}^A=e_{\beta_1}^{(a_1)}e_{\beta_2}^{(a_2)}\cdots e_{\beta_l}^{(a_l)}.
\end{gather}
Then $\{E_{\bf i}^A\mid A\in(\Z_{\ge0})^l\}$ forms a basis of $U_q^+(\geh)$.
We hope that the notations $e_{i_r}$ with $i_r \in I$ and
$e_{\beta_r}$ with a positive root $\beta_r$
can be distinguished properly from the context.
In particular $e^{(a_r)}_{\beta_r}
= (e_{\beta_r})^{a_r}/\prod\limits_{m=1}^{a_r}\frac{p_r^m-p^{-m}_r}{p_r-p_r^{-1}}$
with $p_r = q^{(\beta_r, \beta_r)/2}$.

\section[Quantized algebra of functions $A_q(\geh)$]{Quantized algebra of functions $\boldsymbol{A_q(\geh)}$}\label{a:sec3}

\subsection{Def\/inition}

Following \cite{Kas93}
we give the def\/inition of the quantized algebra of functions $A_q(\geh)$.
It is valid for any symmetrizable Kac--Moody algebra $\geh$. Let $O_\mathrm{int}(\geh)$ be the
category of integrable left $U_q(\geh)$-modules $M$ such that, for any $u\in M$, there exists $l\ge0$
satisfying $e_{i_1}\cdots e_{i_l}u=0$ for any $i_1,\ldots,i_l\in I$. Then $O_\mathrm{int}(\geh)$
is semisimple and any simple object is isomorphic to the irreducible module $V(\la)$ with
dominant integral highest weight $\la$.
Similarly, we can consider the category $O_\mathrm{int}(\geh^\mathrm{opp})$ of integrable
right $U_q(\geh)$-modules $M^r$ such that, for any $v\in M^r$, there exists $l\ge0$ satisfying
$vf_{i_1}\cdots f_{i_l}=0$ for any $i_1,\ldots,i_l\in I$. $O_\mathrm{int}(\geh^\mathrm{opp})$
is also semisimple and any simple object is isomorphic to the irreducible module $V^r(\la)$ with
dominant integral highest weight~$\la$. Let~$u_\la$ (resp.~$v_\la$) be a highest-weight vector
of~$V(\la)$ (resp.~$V^r(\la)$). Then there exists a unique bilinear form $(\,,\,)$
\[
V^r(\la)\ot V(\la)\rightarrow \Q(q)
\]
satisfying
\begin{gather*}
(v_\la,u_\la) =1\qquad\text{and}\qquad
(vP,u) =(v,Pu)\quad\text{for} \ \ v\in V^r(\la), \ \ u\in V(\la), \ \ P\in U_q(\geh).
\end{gather*}

Let $U_q(\geh)^*$ be $\mathrm{Hom}_{\Q(q)}(U_q(\geh),\Q(q))$ and $\pair{\,,\,}$ be the canonical
pairing between $U_q(\geh)^*$ and $U_q(\geh)$. The comultiplication $\Delta$ of $U_q(\geh)$ induces
a multiplication of $U_q(\geh)^*$ by
\begin{gather}\label{eq:mult}
\pair{\varphi\varphi',P}=\pair{\varphi\otimes\varphi',\Delta(P)}\qquad\text{for} \ \ P\in U_q(\geh),
\end{gather}
thereby giving $U_q(\geh)^*$ the structure of $\Q(q)$-algebra. It also has a $U_q(\geh)$-bimodule structure
by
\begin{gather}\label{eq:bimodule}
\pair{x\varphi y,P}=\pair{\varphi,yPx}\qquad\text{for} \ \ x,y,P\in U_q(\geh).
\end{gather}
We def\/ine the subalgebra $A_q(\geh)$ of $U_q(\geh)^*$ by
\[
A_q(\geh)=\big\{\varphi\in U_q(\geh)^*;U_q(\geh)\varphi\text{ belongs to }O_\mathrm{int}(\geh)\text{ and }
\varphi U_q(\geh)\text{ belongs to }O_\mathrm{int}(\geh^\mathrm{opp})\big\},
\]
and call it the quantized algebra of functions.

The following theorem is the $q$-analogue of the Peter--Weyl theorem.
See e.g.~\cite{Kas93}  for a proof.

\begin{theorem}\label{th:peter}
As a $U_q(\geh)$-bimodule $A_q(\geh)$ is isomorphic to $\bigoplus_\la V^r(\la)\ot V(\la)$, where $\la$ runs
over all dominant integral weights, by the homomorphisms
\[
\Psi_\la: \ V^r(\la)\ot V(\la)\rightarrow A_q(\geh)
\]
given by
\[
\langle \Psi_\la(v\ot u),P\rangle = (v,Pu)
\]
for $v\in V^r(\la)$, $u\in V(\la)$, and $P\in U_q(\geh)$.
\end{theorem}

Let us now assume that $\geh$ is a f\/inite-dimensional simple Lie algebra. Then $A_q(\geh)$ turns out
a~Hopf algebra. See e.g.~\cite[Chapter~9]{Joseph}. Its comultiplication is also denoted by~$\Delta$.

Let $\mathcal{R}$ be the
universal $R$ matrix for $U_q(\geh)$. For its explicit formula see e.g.~\cite[p.~273]{CPbook}. For our purpose
it is enough to know that
\begin{gather}\label{eq:univR}
\mathcal{R}\in q^{(\text{wt}\,\cdot,\text{wt}\,\cdot)}
\bigoplus_{\beta\in Q^+}(U_q^+)_\beta\ot(U_q^-)_{-\beta},
\end{gather}
where $q^{(\text{wt}\,\cdot,\text{wt}\,\cdot)}$ is an operator acting on
the tenor product $u_\la\ot u_\mu$ of weight vectors $u_\la$, $u_\mu$ of weight $\la$, $\mu$ by
$q^{(\text{wt}\,\cdot,\text{wt}\,\cdot)}(u_\la\ot u_\mu)=q^{(\la,\mu)}
u_\la\ot u_\mu$, $Q_+=\bigoplus_i\Z_{\ge0}\alpha_i$, and $(U_q^\pm)_{\pm\beta}$ is the subspace of
$U_q^\pm(\geh)$ spanned by root vectors corresponding to $\pm\beta$.

Fix $\la$, let $\{u^\la_j\}$ and $\{v^\la_i\}$ be bases of $V(\la)$ and $V^r(\la)$ such that
$(v^\la_i,u^\la_j)=\delta_{ij}$, and $\varphi^\la_{ij}=\Psi_\la(v^\la_i \otimes u^\la_j)$. Let $R$ be the
so-called constant $R$ matrix for $V(\la)\ot V(\mu)$. Denoting the homomorphism $U_q(\geh)\rightarrow
\mathrm{End}(V(\la))$ by $\pi_\la$, it is given as
\begin{gather}\label{eq:constR}
R\propto(\pi_\la\ot\pi_\mu)(\sigma\mathcal{R}),
\end{gather}
where $\sigma$ stands for the exchange of the f\/irst and second components.
The scalar multiple is determined appropriately depending on~$\geh$. The reason we apply $\sigma$
is that it agrees to the convention of~\cite{RTF}. $R$ satisf\/ies
\[
R\Delta(x)=\Delta'(x)R\qquad\text{for any} \ \ x\in U_q(\geh),
\]
where $\Delta'=\sigma\circ\Delta$.
Def\/ine matrix elements $R_{ij,kl}$ by $R(u^\la_k\ot u^\mu_l)=\sum_{i,j}R_{ij,kl}u^\la_i\ot u^\mu_j$.
Def\/ine the right action of $R$ on $V^r(\la)\ot V^r(\mu)$ in such a way that
$((v^\la_i\ot v^\mu_j)R,u^\la_k\ot u^\mu_l)=(v^\la_i\ot v^\mu_j,R(u^\la_k\ot u^\mu_l))$
holds. Then we have $(v^\la_i\ot v^\mu_j)R=\sum_{k,l}R_{ij,kl}v^\la_k\ot v^\mu_l$. From
\begin{gather*}
\sum_{m,p}R_{ij,mp}\pair{\varphi^\la_{mk}\varphi^\mu_{pl},x}
 = \sum_{m,p}\! R_{ij,mp}\pair{\varphi^\la_{mk}\ot\varphi^\mu_{pl},\Delta(x)}
 =\sum_{m,p} \! R_{ij,mp}(v^\la_m\ot v^\mu_p,\Delta(x)(u^\la_k\ot u^\mu_l))\\
 \hphantom{\sum_{m,p}R_{ij,mp}\pair{\varphi^\la_{mk}\varphi^\mu_{pl},x}}{}
 =((v^\la_i\ot v^\mu_j)R,\Delta(x)(u^\la_k\ot u^\mu_l))
 =(v^\la_i\ot v^\mu_j,R\Delta(x)(u^\la_k\ot u^\mu_l))\\
 \hphantom{\sum_{m,p}R_{ij,mp}\pair{\varphi^\la_{mk}\varphi^\mu_{pl},x}}{}
 =\sum_{m,p}(v^\la_i\ot v^\mu_j,\Delta'(x)(u^\la_m\ot u^\mu_p))R_{mp,kl}\\
\hphantom{\sum_{m,p}R_{ij,mp}\pair{\varphi^\la_{mk}\varphi^\mu_{pl},x}}{}
 =\sum_{m,p}(v^\mu_j\ot v^\la_i,\Delta(x)(u^\mu_p\ot u^\la_m))R_{mp,kl}\\
 \hphantom{\sum_{m,p}R_{ij,mp}\pair{\varphi^\la_{mk}\varphi^\mu_{pl},x}}{}
 =\sum_{m,p}\pair{\varphi^\mu_{jp}\ot\varphi^\la_{im},\Delta(x)}R_{mp,kl}
=\sum_{m,p}\pair{\varphi^\mu_{jp}\varphi^\la_{im},x}R_{mp,kl}
\end{gather*}
for any $x \in U_q(\geh)$, we have
\begin{gather}\label{eq:RTT}
\sum_{m,p}R_{ij,mp}\varphi^\la_{mk}\varphi^\mu_{pl}=\sum_{m,p}\varphi^\mu_{jp}\varphi^\la_{im}R_{mp,kl}.
\end{gather}
We call such a relation $RTT$ relation.

\subsection[Right quotient ring $A_q(\geh)_\mathcal{S}$]{Right quotient ring $\boldsymbol{A_q(\geh)_\mathcal{S}}$}

For later use we require a certain right quotient ring of $A_q(\geh)$ by a suitable multiplicatively
closed subset~$\mathcal{S}$. We f\/irst review the general construction from \cite[Chapter~2]{McCoRob}.

Let $R$ be a noncommutative ring with $1$ and $\mathcal{S}$ a multiplicatively closed subset of~$R$.
The following condition is called the right Ore condition:
\begin{itemize}\itemsep=0pt
\item[(Ore)] For any $r\in R$, $s\in\mathcal{S}$, $r\mathcal{S}\cap sR\neq\varnothing$.
\end{itemize}
Set
\[
\mathrm{ass}\,\mathcal{S}=\{r\in R\mid rs=0\mbox{ for some }s\in\mathcal{S}\}.
\]
Then under the right Ore condition $\mathrm{ass}\,\mathcal{S}$ turns out a two-sided ideal.
Let $\overline{\phantom{a}}:R\rightarrow R/\mathrm{ass}\,\mathcal{S}$ denote the canonical projection.
Suppose
\begin{itemize}\itemsep=0pt
\item[(reg)] $\overline{\mathcal{S}}$ consists of regular elements, namely, elements $x$ such that
both $xr=0$ and $rx=0$ imply $r=0$.
\end{itemize}
Then a theorem in \cite[Chapter~2]{McCoRob} states

\begin{theorem}[Theorem~2.1.12 of~\cite{McCoRob}] \label{th:McCoRob}
The right quotient ring $R_\mathcal{S}$ exists, if and only if $($Ore$)$ and $($reg$)$ are satisfied.
\end{theorem}

By passing to the images by $\overline{\phantom{a}}$, it suf\/f\/ices to consider the case when
$\mathrm{ass}\,\mathcal{S}=0$, and then elements of $R_\mathcal{S}$ are of the form~$r/s$.
For $r_i/s_i\in R/\mathcal{S}$ ($i=1,2$) the addition and multiplication formulae are given by
\begin{gather}\label{add-mult}
r_1/s_1+r_2/s_2=(r_1u+r_2u')/(s_1u),\qquad (r_1/s_1)(r_2/s_2)=(r_1v')/(s_2v),
\end{gather}
where $u$, $u'$, $v$, $v'$ are so chosen that $s_1u=s_2u'$ ($u\in\mathcal{S},u'\in R$), $r_2v=s_1v'$
($v\in\mathcal{S},v'\in R$).

Let us return to our case where $R=A_q(\geh)$.
\begin{definition}\label{def:sigma}
For any $i \in I$, let $u_{w_0\varpi_i}$ (resp.~$v_{\varpi_i}$) be a lowest (resp. highest) weight vector
of $V(\varpi_i)$ (resp. $V^r(\varpi_i)$). Set
\[
\sigma_i=\Psi_{\varpi_i}(v_{\varpi_i}\ot u_{w_0\varpi_i}).
\]
\end{definition}

The following proposition is proven in~\cite{Joseph}. However, we dare to prove again, since conventions
might be dif\/ferent.

\begin{proposition}[Corollary~9.1.4 of \cite{Joseph}] \label{pr:q-comm}
Let $\varphi_{\la\mu}$ be an element of $A_q(\geh)$ such that $k_i\varphi_{\la\mu}=
q_i^{\pair{h_i,\mu}}\varphi_{\la\mu},\varphi_{\la\mu}k_i=q_i^{\pair{h_i,\la}}\varphi_{\la\mu}$
for any $i \in I$. Then the following commutation relation holds:
\[
q^{(\varpi_i,\la)}\sigma_i\varphi_{\la\mu}=q^{(w_0\varpi_i,\mu)}\varphi_{\la\mu}\sigma_i.
\]
In particular, $\sigma_i\sigma_j=\sigma_j\sigma_i$ for any $i$, $j$.
\end{proposition}

\begin{proof}
Without loss of generality one can assume $\varphi_{\la\mu}=\Psi_\nu(v_\la\ot u_\mu)$ for some
$\nu,v_\la\in V^r(\nu)$, $u_\mu\in V(\nu)$ such that $k_iu_\mu=q_i^{\pair{h_i,\mu}}u_\mu$,
$v_\la k_i=q_i^{\pair{h_i,\la}}v_\la$. In view of \eqref{eq:univR}, \eqref{eq:constR} we have
\[
R(u_{w_0\varpi_i}\ot u_\mu)=q^{(w_0\varpi_i,\mu)}u_{w_0\varpi_i}\ot u_\mu,\qquad
(v_{\varpi_i}\ot v_\la)R=q^{(\varpi_i,\la)}v_{\varpi_i}\ot v_\la.
\]
Then \eqref{eq:RTT} implies the commutation relation.
The second relation follows from the f\/irst one, since $(\varpi_i,\varpi_i)
=(w_0\varpi_i,w_0\varpi_i)$.
\end{proof}

Let $n$ be the rank of $\geh$ and def\/ine
\begin{gather*}
\mathcal{S}=\big\{\sigma_1^{m_1}\cdots\sigma_n^{m_n}\mid m_1,\ldots,m_n\in\Z_{\ge0}\big\},
\end{gather*}
which is obviously multiplicatively closed subset of $A_q(\geh)$.

\begin{lemma}\label{le:ff}
Let $s$ be a nonzero element in $\mathrm{Im}\,\Psi_\la$
satisfying $f_is=sf_i=0$ for any $i \in I$.
Then $s\in\Q(q)^\times\sigma_1^{\la_1}\cdots\sigma_n^{\la_n}$
where $\Q(q)^\times=\Q(q)\setminus\{0\}$ and
$\la_i=\pair{h_i,\la}$.
\end{lemma}
\begin{proof}
By \eqref{eq:mult}, \eqref{eq:bimodule} $f_i\sigma_1^{\la_1}\cdots\sigma_n^{\la_n}
=\sigma_1^{\la_1}\cdots\sigma_n^{\la_n}f_i=0$ for any $i$ and $\sigma_1^{\la_1}\cdots\sigma_n^{\la_n}$
belongs to $\mathrm{Im}\,\Psi_\la$. By Theorem~\ref{th:peter} such an element is unique up to an element
of $\Q(q)^\times$.
\end{proof}

In particular
$\sigma_i$ is characterized as the unique element (up to an overall constant)
in $\mathrm{Im}\, \Psi_{\varpi_i}$ such that
$f_j \sigma_i = \sigma_i f_j = 0$ for all $j \in I$.
We remark that Theorem~\ref{th:peter} implies that
if a nonzero element $\varphi_{\lambda,\mu} \in A_q(\geh)$
satisf\/ies the assumption of Proposition~\ref{pr:q-comm} and
$f_j \varphi_{\lambda,\mu}
= \varphi_{\lambda,\mu} f_j = 0$ for all~$j$,
then $\lambda = w_0\mu$ must hold.

In \cite{Joseph} it is shown that $A_q(\geh)$ is an integral domain~(Lemma 9.1.9), hence (reg) is
satisf\/ied, and that~(Ore) is also satisf\/ied (Lemma~9.1.10). Therefore we have the following theorem.
(A~proof is attached for self-containedness.)

\begin{theorem}\label{th:quotient}
The right quotient ring $A_q(\geh)_\mathcal{S}$ exists.
\end{theorem}

\begin{proof}
In view of Theorem \ref{th:McCoRob} it is enough to show that
\begin{itemize}\itemsep=0pt
\item[(1)] if $\varphi\ne0$, then $\varphi s\ne0$ for any $s\in\mathcal{S}$,
\item[$(1')$] if $\varphi\ne0$, then $s\varphi\ne0$ for any $s\in\mathcal{S}$, and
\item[(2)] the right Ore condition is satisf\/ied,
\end{itemize}
since (1) implies $\mathrm{ass}\,\mathcal{S}=0$, then~(1) and~$(1')$ imply $\overline{\mathcal{S}}=
\mathcal{S}$ consists of regular elements.

Let us prove~(1). Let $\varphi=\sum_j\varphi_j$ be the two-sided weight decomposition.
If $\varphi_js\ne0$ for some~$j$,
$\varphi s\ne0$ since the weights of $\varphi_js$ are distinct.
Hence we can reduce the claim when $\varphi$ is a~weight vector.
Suppose $\varphi=\sum_\mu\varphi_\mu$, $\varphi_\mu\in\mathrm{Im}\,\Psi_\mu$ and let $\la$ be a maximal
weight, with respect to the standard ordering on weights, such that $\varphi_\la\ne0$. Choose sequences
$i_1,\ldots,i_k$ and $j_1,\ldots,j_l$ such that
$f_{i_k}\cdots f_{i_1}\varphi_\la f_{j_1}\cdots f_{j_l}$
turns out a left-lowest and right-highest weight vector.
Then by Lemma~\ref{le:ff} it coincides with $cs'$
with some $c\in\Q(q)^\times$, $s'\in\mathcal{S}$. Then
\[
f_{i_k}\cdots f_{i_1}(\varphi s)f_{j_1}\cdots f_{j_l}=c's's+\cdots
\]
with another $c'\in\Q(q)^\times$. By the maximality of $\la$ the remaining part $+\cdots$ in the right-hand
side does not contain
the terms with the same two-sided weight. Hence $\cdots=0$.
Therefore, the left-hand side is
not $0$ and we conclude $\varphi s\ne0$.

$(1')$ is similar.
For (2) we can reduce the claim when $\varphi$ is a weight vector, and in this case the claim is clear from
Proposition~\ref{pr:q-comm}.
\end{proof}

\subsection{Realization by generators and relations}

We consider the fundamental representation $V(\varpi_1)$ of $U_q(\geh)$ for $\geh=A_{n-1},C_n,G_2$.
Set $N=\dim V(\varpi_1)$. It is known \cite{D86,RTF} that $A_q(\geh)$ for $\geh=A_{n-1},C_n,G_2$ is
realized as an associative algebra with appropriate generators $(t_{ij})_{1\le i,j \le N}$
corresponding to $V^r(\varpi_1)\ot V(\varpi_1)$ satisfying $RTT$ relations
\begin{gather}\label{a:rtt}
 \sum_{m,p}R_{ij,mp}t_{mk}t_{pl}
= \sum_{m,p}t_{jp}t_{im}R_{mp,kl},
\end{gather}
and additional ones depending on $\geh$.
See below for each $\geh$ under consideration.
In all cases, there exists a comultiplication $\Delta: A_q \rightarrow A_q \otimes A_q$ given by
\begin{gather}\label{a:cpro}
\Delta(t_{ij}) = \sum_k t_{ik}\otimes t_{kj}.
\end{gather}

\subsubsection[$A_{n-1}$ case]{$\boldsymbol{A_{n-1}}$ case}\label{a:sss:a1}

We present formulae for $A_q(A_{n-1})$. In this case $N=n$.
Let $u_1$ and $v_1$ be the highest-weight vectors of
$V(\varpi_1)$ and $V^r(\varpi_1)$ such that $(v_1,u_1)=1$ and
set $u_j=f_{j-1}f_{j-2}\cdots f_1u_1$, $v_j=v_1e_1e_2\cdots e_{j-1}$
for $2\le j\le n$. Then the constant $R$ matrix is given by
\begin{gather*}
\sum_{i,j,k,l}R_{ij,kl}E_{ik}\otimes E_{jl}=
q\sum_i E_{ii}\otimes E_{ii} + \sum_{i\neq j}E_{ii}\otimes E_{jj}
+\big(q-q^{-1}\big)\sum_{i>j}E_{ij}\otimes E_{ji},
\end{gather*}
where $E_{ij}$ is the matrix unit. Def\/ine $t_{ij}=\Psi_{\varpi_1}(v_i \otimes u_j)$.
Then the $RTT$ relations among $(t_{ij})_{1\le i,j\le N}$ read explicitly as follows
\begin{gather*}
 [t_{ik}, t_{jl}]=\begin{cases}0, & i<j, \ k>l,\\
\big(q-q^{-1}\big)t_{jk}t_{il}, &  i<j, \ k<l,
\end{cases}\\
 t_{ik}t_{jk} = q t_{jk}t_{ik},\quad i<j,\qquad
t_{ki}t_{kj} = q t_{kj}t_{ki},\quad i<j.
\end{gather*}
In $A_{n-1}$ case we need another condition that the quantum determinant is 1, i.e.,
\begin{gather*}
\sum_{\sigma \in \mathfrak{S}_n}(-q)^{\ell(\sigma)}
t_{1\sigma_1}\cdots t_{n\sigma_n} = 1,
\end{gather*}
where $\mathfrak{S}_n=W(A_{n-1})$ is the symmetric group of degree $n$ and
$\ell(\sigma)$ is the length of $\sigma$.

According to Def\/inition \ref{def:sigma}, we have
$\sigma_1 = t_{13}$ and $\sigma_2 = t_{12}t_{23}-qt_{22}t_{13}$.
As an exposition,  we note that
$\sigma_i e_i$ in (\ref{a:kobetsu}) is  derived from
\begin{gather*}
\langle \sigma_1e_1, P\rangle  =\langle t_{13}e_1, P\rangle
= (v_1e_1, Pu_3) = (v_2, Pu_3) = \langle t_{23}, P \rangle,\\
\langle \sigma_2e_2, P\rangle
= \langle (t_{12}\otimes t_{23}-qt_{22}\otimes t_{13})\Delta(e_2), \Delta(P)\rangle
=\langle t_{12}k_2 \otimes t_{23}e_2- q t_{22}e_2\otimes t_{13}, \Delta(P)\rangle\\
\hphantom{\langle \sigma_2e_2, P\rangle}{}
= \langle t_{12} \otimes t_{33}- q t_{32}\otimes t_{13}, \Delta(P)\rangle
= \langle t_{12}t_{33}- q t_{32} t_{13}, P\rangle
\end{gather*}
for any $P \in U_q(A_2)$.
See e.g.~\cite{NYM} for an extensive treatment.

\subsubsection[$C_n$ case]{$\boldsymbol{C_n}$ case}\label{a:sss:c1}

We present formulae for $A_q(C_n)$. In this case $N=2n$.
Let $u_1$ be the highest-weight vector of~$V(\varpi_1)$ and
def\/ine $u_j$ for $2\le j\le 2n$ recursively by $u_{j+1}=f_ju_j$~$(j\le n)$, $-f_{2n-j}u_j$~$(j>n)$.
Let~$\{v_i\}$ be the dual basis to $\{u_i\}$ in $V^r(\varpi_1)$,
namely, $\{v_i\}$ are determined
by $(v_i,u_j)=\delta_{ij}$. Then the constant $R$ matrix is given by
\begin{gather*}
 \sum_{i,j,k,l}R_{ij,kl}E_{ik}\otimes E_{jl}=
q\sum_i E_{ii}\otimes E_{ii} + \sum_{i\neq j, j'}E_{ii}\otimes E_{jj}
+q^{-1}\sum_i E_{ii}\otimes E_{i'i'}\\
\hphantom{\sum_{i,j,k,l}R_{ij,kl}E_{ik}\otimes E_{jl}=}{}
+\big(q-q^{-1}\big)\sum_{i>j}E_{ij}\otimes E_{ji}
-\big(q-q^{-1}\big)\sum_{i>j}\epsilon_i\epsilon_jq^{\varrho_i-\varrho_j}E_{ij}\otimes E_{i'j'},\\
i' = 2n+1-i,
\qquad
\epsilon_i = 1, \quad 1\le i \le n,\qquad \epsilon_i = -1, \quad n<i \le 2n,\\
(\varrho_1,\ldots, \varrho_{2n}) = (n-1,n-2,\ldots,1,0,0,-1,\ldots, -n+1).
\end{gather*}
Def\/ine $t_{ij}=\Psi_{\varpi_1}(v_i \otimes u_j)$. The $RTT$ relations are given by \eqref{a:rtt}
with the above $R_{ij,kl}$. Additional relations are given by
\begin{gather*}
 \sum_{j,k,l}C_{jk}C_{lm}t_{ij}t_{lk}
= \sum_{j,k,l}C_{ij}C_{kl}t_{kj}t_{lm} = -\delta_{im},
\qquad  C_{ij}= \delta_{i,j'}\epsilon_i q^{\varrho_j}.
\end{gather*}

\subsubsection[$G_2$ case]{$\boldsymbol{G_2}$ case}\label{a:sss:g1}

We have $N=7$ in this case. We adopt the basis $\{u_i\}$ of $V(\varpi_1)$ that has the
representation matrices given as in \cite[equation~(29)]{Sasaki}, and let $\{v_i\}$ the dual basis
in $V^r(\varpi_1)$. Def\/ine $t_{ij}=\Psi_{\varpi_1}(v_i \otimes u_j)$. Then
$A_q(G_2)$ is generated by $(t_{ij})_{1\le i,j \le 7}$ satisfying~(i) and~(ii) given below.
\begin{itemize}\itemsep=0pt
\item[(i)] $RTT$ relations (\ref{a:rtt})
with the structure constants specif\/ied by
$R_{ij,kl} = R^{ij}_{\;\;kl}$ in \cite[equation~(33)]{Sasaki}.
\item[(ii)] Additional relations
\begin{gather}\label{a:re4}
g^{ij} = \sum_{k,l}t_{jl}t_{ik}g^{kl},\qquad
\sum_k f^{ij}_{\;\;k} t_{km} = \sum_{k,l}t_{jl}t_{ik} f^{k l}_{\;\;m},
\end{gather}
where $g^{ij}$ and $f^{ij}_{\;\;k}$ are
given by \cite[equations~(30), (31)]{Sasaki}.
\end{itemize}
The relations \cite[equations~(20), (22)]{Sasaki} are
equivalent to~(\ref{a:re4})  if the~$RTT$ relations are imposed.
See the explanation after \cite[Def\/inition~7]{Sasaki}.
Note also that we use the opposite indices of the Dynkin diagram to~\cite{Sasaki}.

\section[Representations of $A_q(\geh)$]{Representations of $\boldsymbol{A_q(\geh)}$}\label{a:sec4}

\subsection{General remarks}\label{a:ss:gr}

Let us recall the results in \cite{So2, VS}
on the representations of $A_q(\mathfrak{g})$
necessary in this paper.
Consider the simplest example $A_q(A_1)$ generated by
$t_{11}$, $t_{12}$, $t_{21}$, $t_{22}$
with the relations
\begin{gather*}
t_{11}t_{21} = qt_{21}t_{11},\qquad
t_{12}t_{22} = qt_{22}t_{12},\qquad
t_{11}t_{12} = qt_{12}t_{11},\qquad
t_{21}t_{22} = qt_{22}t_{21},\\
[t_{12},t_{21}]=0,\quad [t_{11},t_{22}]=(q-q^{-1})t_{21}t_{12},
\qquad t_{11}t_{22}-qt_{12}t_{21}=1.
\end{gather*}
Let $\mathrm{Osc}_q =
\langle  {\rm{\bf a}}^+, {\rm{\bf a}}^-, {\rm{\bf k}} \rangle$
be the $q$-oscillator algebra, i.e.,
an associative algebra with the relations
\begin{gather}\label{a:qoa}
{\rm{\bf k}}  {\rm{\bf a}}^+ = q {\rm{\bf a}}^+{\rm{\bf k}},\qquad
{\rm{\bf k}}\,{\rm{\bf a}}^- = q^{-1}{\rm{\bf a}}^-{\rm{\bf k}},\qquad
{\rm{\bf a}}^- {\rm{\bf a}}^+ = {\bf 1}-q^2{\rm{\bf k}}^2,\qquad
{\rm{\bf a}}^+{\rm{\bf a}}^- = {\bf 1}-{\rm{\bf k}}^2.
\end{gather}
It has a representation on
the Fock space
$\mathcal{F}_q = \bigoplus_{m\ge 0}{\mathbb C}(q)|m\rangle$:
\begin{gather}\label{a:akf}
{\rm{\bf k}}|m\rangle = q^m |m\rangle,\qquad
{\rm{\bf a}}^+|m\rangle = |m+1\rangle,\qquad
{\rm{\bf a}}^-|m\rangle = (1-q^{2m})|m-1\rangle.
\end{gather}
In what follows, the symbols
${\rm{\bf k}}$, ${\rm{\bf a}}^+$, ${\rm{\bf a}}^-$ shall
also be regarded as the elements from $\mathrm{End}(\mathcal{F}_q)$.
It is easy to check that the following map $\pi$
def\/ines an irreducible representation of
$A_q(A_1)$ on $\mathcal{F}_q$:
\begin{gather}\label{a:pisl2}
\pi:  \ \begin{pmatrix}
t_{11} & t_{12}\\
t_{21} & t_{22}
\end{pmatrix} \mapsto
\begin{pmatrix}
\mu\,{\rm {\bf a}}^- &  \alpha {\rm {\bf k}}\\
-q\alpha^{-1}{\rm {\bf k}} & \mu^{-1}{\rm {\bf a}}^+
\end{pmatrix},
\end{gather}
where $\alpha$, $\mu$ are nonzero parameters.

\begin{theorem}[\cite{So2,VS}]\label{a:th:so}
\qquad
\begin{enumerate}\itemsep=0pt
\item[$(1)$]  For each vertex $i$
of the Dynkin diagram of $\mathfrak{g}$,
$A_q(\mathfrak{g})$ has an irreducible
representa\-tion~$\pi_i$ factoring through
\eqref{a:pisl2} via
$A_q(\mathfrak{g}) \twoheadrightarrow
A_{q_i}(sl_{2,i})$.
$(sl_{2,i}$ denotes the $sl_2$-subalgebra of $\mathfrak{g}$
asso\-ciated to~$i.)$

\item[$(2)$] Irreducible representations of $A_q(\mathfrak{g})$ are
in one to one correspondence with the elements of the Weyl group
$W$ of $\mathfrak{g}$.

\item[$(3)$] Let $w =s_{i_1}\cdots s_{i_l} \in W$ be an
reduced expression in terms of the simple reflections.
Then the irreducible representation corresponding to $w$ is
isomorphic to $\pi_{i_1}\otimes \cdots \otimes \pi_{i_l}$.
\end{enumerate}
\end{theorem}

Actually the assertions (2) and (3) hold
up to the degrees of freedom of the parame\-ters~$\alpha$,~$\mu$ in~(\ref{a:pisl2}).
See~\cite{So2} for the detail.
We call $\pi_i$ $(i = 1,\ldots, \mathrm{rank}\, \mathfrak{g})$
the {\em fundamental representations}.
For simplicity we denote
$\pi_{i_1} \otimes \cdots \otimes \pi_{i_l}$
by $\pi_{i_1,\ldots, i_l}$.

A crucial corollary of Theorem \ref{a:th:so} is the following:
\begin{gather*}
\mathrm{If } \ s_{i_1}\cdots s_{i_l} = s_{j_1}\cdots s_{j_l}
\in W \
\text{are reduced expressions, then } \
\pi_{i_1, \ldots, i_l} \simeq \pi_{j_1, \ldots, j_l}.
\end{gather*}
In particular, there exists the isomorphism
$\Phi: \mathcal{F}_{q_{i_1}}\otimes \cdots \otimes \mathcal{F}_{q_{i_l}}
\rightarrow
\mathcal{F}_{q_{j_1}}\otimes \cdots \otimes \mathcal{F}_{q_{j_l}}$
characterized (up to an overall constant) by
\begin{gather*}
\pi_{j_1, \ldots, j_l}(g)\circ \Phi = \Phi \circ \pi_{i_1, \ldots, i_l}(g)\qquad
\forall\, g \in A_q(\mathfrak{g}).
\end{gather*}
Here $\pi_{i_1, \ldots, i_l}(g=t_{ij})$ for example
means the tensor product representation
$\sum_{r_1,\ldots, r_{l-1}}
\pi_{i_1}(t_{i r_1}) \otimes \cdots \otimes \pi_{i_l}(t_{r_{l-1},j})$
obtained by the $(l-1)$-fold application of the
coproduct~(\ref{a:cpro}).

Elements of the Fock space
$|m_1\rangle \otimes \cdots \otimes |m_l\rangle
\in \mathcal{F}_{q_{j_1}}\otimes \cdots \otimes \mathcal{F}_{q_{j_l}}$
will simply be denoted by
$|m_1, \ldots, m_l\rangle$.
We will always normalize the intertwiner by the condition
$\Phi |0,0,\ldots, 0\rangle =  |0,0,\ldots, 0\rangle$.
The exchange of the $i$th and the $j$th tensor components from the left
will be denoted by~$P_{ij}$.
In the remainder of this section we concentrate on
$A_q(\mathfrak{g})$ of rank~$2$ cases
$\mathfrak{g}=A_2, C_2$ and $G_2$,
and present the concrete forms of the
fundamental representations,
def\/inition of the intertwiners with
a few examples of their matrix elements.

\subsection[$A_2$ case]{$\boldsymbol{A_2}$ case}\label{a:ss:a2}

Let $T=(t_{ij})_{1\le i,j \le 3}$ be the
$3\times 3$ matrix of the generators of $A_q(A_2)$.
The fundamental representations
$\pi_i: A_q(A_2) \rightarrow \mathrm{End}(\mathcal{F}_q)$
$(i=1,2)$ are given by
\begin{gather}\label{a:piTa}
\pi_1(T) =
\begin{pmatrix}
\mu_1{\rm {\bf a}}^- &  \alpha_1{\rm {\bf k}} & 0 \\
-q\alpha^{-1}_1{\rm {\bf k}} & \mu^{-1}_1{\rm {\bf a}}^+ & 0\\
0 & 0 & 1
\end{pmatrix},
\qquad
\pi_2(T) =
\begin{pmatrix}
1 & 0 & 0\\
0 & \mu_2{\rm {\bf a}}^- &  \alpha_2{\rm {\bf k}}\\
0 & -q\alpha^{-1}_2{\rm {\bf k}} & \mu^{-1}_2{\rm {\bf a}}^+
\end{pmatrix},
\end{gather}
where $\alpha_i$, $\mu_i$ are nonzero parameters.

The Weyl group $W=\langle s_1, s_2 \rangle$ is the Coxeter system
with the relations
\begin{gather*}
s_1^2=s_2^2=1,\qquad
s_1s_2s_1 = s_2s_1s_2.
\end{gather*}
Thus we have the isomorphism
$\pi_{121} \simeq \pi_{212}$.
Let $\Phi$ be the corresponding intertwiner and denote by $R$
the checked intertwiner $\Phi^\vee$ explained after (\ref{a:main})
\begin{gather*}
 \pi_{121} \Phi = \Phi  \pi_{212},\qquad
\pi_{121} R = R  \pi'_{212},
\qquad \pi'_{212} = P_{13}\pi_{212}P_{13},\qquad
 R = \Phi P_{13}\in \mathrm{End}\big(\mathcal{F}_q^{\otimes 3}\big).
\end{gather*}
For example
$\pi'_{212}(t_{ij}) = \sum_{k,l}
\pi_2(t_{l,j})\otimes \pi_1(t_{k,l})\otimes \pi_2(t_{ik})$.
Def\/ine the matrix elements of $R$ and its parameter-free part
$\mathscr{R}$ by
\begin{gather*}
R |i,j,k \rangle=
\sum_{a,b,c}R^{abc}_{ijk}|a,b,c \rangle,
\qquad
R^{abc}_{ijk}  = \mu_1^{a-j+k}\mu_2^{b-a-k}
\mathscr{R}^{abc}_{ijk}.
\end{gather*}
Then the following properties are valid for
$\mathscr{R} = (\mathscr{R}^{abc}_{ijk})$~\cite{KO1}:
\begin{gather}
 \mathscr{R}^{abc}_{ijk} \in \mathbb{Z}[q],\qquad
\mathscr{R}^{abc}_{ijk}=0 \quad \mathrm{unless}\ \ (a+b, b+c)=(i+j,j+k),
\label{a:claw}\\
 \mathscr{R}^{-1} = \mathscr{R},
\qquad\mathscr{R}^{abc}_{ijk} = \mathscr{R}^{cba}_{kji},\qquad
{\mathscr{R}}^{abc}_{ijk}
= \frac{(q^2)_i(q^2)_j(q^2)_k}{(q^2)_a(q^2)_b(q^2)_c}
{\mathscr{R}}^{ijk}_{abc},\label{a:rprop}\\
\mathscr{R}^{abc}_{ijk}|_{q=0} = \delta_{i, b+(a-c)_+}\delta_{j,\min(a,c)}
\delta_{k, b+(c-a)_+}.\label{a:r0}
\end{gather}
Here $(q^2)_a = \prod\limits_{m=1}^a(1-q^{2m})$
and $(y)_+ = \max(0,y)$.
Due to~(\ref{a:claw}), $\mathscr{R}$ is the inf\/inite direct sum
of f\/inite-dimensional matrices.
An explicit formula of $\mathscr{R}^{abc}_{ijk}$ was obtained in~\cite{KV} (unfortunately with misprint) and in
\cite[equation~(59)]{BMS} (in a dif\/ferent context and gauge including square roots).
The formula exactly matching the present convention is
\cite[equation~(2.20)]{KO1}.
The $\mathscr{R}$ satisf\/ies~\cite{KV}
the tetrahedron equation~(\ref{a:teq}).

\begin{example}\label{a:ex:a21}
The following is the list of all the nonzero $\mathscr{R}^{abc}_{314}$:
\begin{gather*}
\mathscr{R}^{041}_{314}  = -q^2 \big(1 - q^4\big) \big(1 - q^6\big) \big(1 - q^8\big), \\
\mathscr{R}^{132}_{314}  =\big(1 - q^6\big) \big(1 - q^8\big)\big(1 - q^4 - q^6 - q^8 - q^{10}\big),\\
\mathscr{R}^{223}_{314}  = q^2 \big(1 + q^2\big) \big(1 + q^4\big) \big(1 - q^6\big) \big(1 - q^6 - q^{10}\big), \\
\mathscr{R}^{314}_{314}  =q^6\big(1 + q^2 + q^4 - q^8 - q^{10} - q^{12} - q^{14}\big),\\
\mathscr{R}^{405}_{314}  =q^{12}.
\end{gather*}
Thus $\mathscr{R}^{abc}_{314}|_{q=0}=\delta_{a,1}\delta_{b,3}\delta_{c,2}$
in agreement with (\ref{a:r0}).
\end{example}

\subsection[$C_2$ case]{$\boldsymbol{C_2}$ case}\label{a:ss:c2}

We have $(q_1,q_2)=(q,q^2)$.
Let $T=(t_{ij})_{1\le i,j \le 4}$ be the
$4\times 4$ matrix of the generators of~$A_q(C_2)$.
We  use
$\mathrm{Osc}_{q^2} =
\langle  {\rm{\bf A}}^+, {\rm{\bf A}}^-, {\rm{\bf K}} \rangle$
in addition to
$\mathrm{Osc}_q =
\langle  {\rm{\bf a}}^+, {\rm{\bf a}}^-, {\rm{\bf k}} \rangle$
(\ref{a:qoa}).
The fundamental representations
$\pi_i: A_q(C_2) \rightarrow \mathrm{End}(\mathcal{F}_{q_i})$
$(i=1,2)$ are given by
\begin{gather}
\pi_1(T) =
\begin{pmatrix}
\mu_1{\rm{\bf a}}^- & \alpha_1{\rm {\bf k}} &  0& 0\\
-q\alpha^{-1}_1{\rm {\bf k}}& \mu^{-1}_1{\rm{\bf a}}^+ &  0& 0\\
0 & 0 & \mu_1{\rm{\bf a}}^- & -\alpha_1{\rm {\bf k}} \\
0 & 0 & q\alpha_1^{-1}{\rm {\bf k}} & \mu_1^{-1}{\rm{\bf a}}^+
\end{pmatrix},
\nonumber\\
\pi_2(T) =
\begin{pmatrix}
1  & 0 & 0 & 0 \\
0 & \mu_2{\rm{\bf A}}^- & \alpha_2{\rm{\bf K}} & 0 \\
0 & -q^2\alpha^{-1}_2{\rm{\bf K}} & \mu_2^{-1}{\rm{\bf A}}^+ & 0 \\
0 & 0 & 0  & 1
\end{pmatrix},\label{a:piTc}
\end{gather}
where $\alpha_i$, $\mu_i$ are nonzero parameters.

The Weyl group $W=\langle s_1, s_2 \rangle$ is the Coxeter system
with the relations
\begin{gather*}
s_1^2=s_2^2=1,\qquad
s_2s_1s_2s_1 = s_1s_2s_1s_2.
\end{gather*}
Thus we have the isomorphism
$\pi_{2121} \simeq \pi_{1212}$.
Let $\Phi$ be the corresponding intertwiner and denote by~$K$
the checked intertwiner $\Phi^\vee$
\begin{gather*}
\pi_{2121} \Phi = \Phi  \pi_{1212},\qquad
\pi_{2121} K = K  \pi'_{2121},\qquad
\pi'_{2121} = P_{14}P_{23}\pi_{1212}P_{14}P_{23},\\
K = \Phi P_{14}P_{23} \in
\mathrm{End}(\mathcal{F}_{q^2}\otimes
\mathcal{F}_{q}\otimes
\mathcal{F}_{q^2}\otimes
\mathcal{F}_{q}).
\end{gather*}
Def\/ine the matrix elements of $K$ and its parameter-free part
$\mathscr{K}$ by
\begin{gather*}
K  |i,j,k,l \rangle=
\sum_{a,b,c,d}K^{abcd}_{ijkl}|a,b,c,d \rangle,
\qquad
K^{abcd}_{ijkl}  = \mu_1^{2(c-k)}\mu_2^{b-j}
\mathscr{K}^{abcd}_{ijkl}.
\end{gather*}
Then the following properties are valid for
$\mathscr{K} = (\mathscr{K}^{abcd}_{ijkl})$ \cite{KO1}:
\begin{gather}
\mathscr{K}^{abcd}_{ijkl} \in \mathbb{Z}[q],\qquad
\mathscr{K}^{abcd}_{ijkl}=0 \quad \mathrm{unless}\ \
(a+b+c, b+2c+d)=(i+j+k, j+2k+l),
\label{a:claw2}\\
\mathscr{K}^{-1} = \mathscr{K},\qquad
\mathscr{K}^{a b c d}_{i j k l}
= \frac{(q^4)_i(q^2)_j(q^4)_k(q^2)_l }{(q^4)_a(q^2)_b(q^4)_c(q^2)_d}
\mathscr{K}^{i j k l}_{a b c d},\label{a:kprop}\\
\mathscr{K}^{abcd}_{ijkl}|_{q=0}
= \delta_{i,a'}\delta_{j,b'}\delta_{k,c'}\delta_{l,d'},\label{a:k0}\\
a'=x + a + b - d,\qquad
b'=c + d - x-\min(a,c + x),\nonumber\\
  c'=\min(a,c+x),\qquad
d'= b+(c-a+x)_+,\qquad x=(c-a+(d-b)_+)_+.\nonumber
\end{gather}
Due to (\ref{a:claw2}), $\mathscr{K}$ is the inf\/inite direct sum
of f\/inite-dimensional matrices.
An explicit formula of $\mathscr{K}^{abcd}_{ijkl}$ is available
in \cite[equations~(3.27), (3.28)]{KO1}.
This $\mathscr{K}$ and $\mathscr{R}$ in Section~\ref{a:ss:a2} satisfy~\cite{KO1}
the 3D ref\/lection equation~(\ref{a:3dref}).

\begin{example}\label{ex:K}
The following is the list of all the nonzero $\mathscr{K}^{abcd}_{2110}$:
\begin{gather*}
\mathscr{K}^{1300}_{2110}=q^8 \big(1 - q^8\big),\\
\mathscr{K}^{2110}_{2110} =-q^4 \big(1 - q^8 + q^{14}\big),\\
\mathscr{K}^{2201}_{2110} =-q^6 \big(1 + q^2\big) \big(1 - q^2 + q^4 - q^6 - q^{10}\big),\\
\mathscr{K}^{3011}_{2110} =1 - q^8 + q^{14},\\
\mathscr{K}^{3102}_{2110} =-q^{10} \big(1 - q + q^2\big) \big(1 + q + q^2\big),\\
\mathscr{K}^{4003}_{2110}=q^4.
\end{gather*}
Thus $\mathscr{K}^{abcd}_{2110}|_{q=0}
= \delta_{a,3}\delta_{b,0}\delta_{c,1}\delta_{d,1}$ in agreement with
(\ref{a:k0}).
\end{example}

\subsection[$G_2$ case]{$\boldsymbol{G_2}$ case}\label{a:ss:g2}

We have $(q_1,q_2)=(q,q^3)$.
Let $T=(t_{ij})_{1\le i,j \le 7}$ be the
$7\times 7$ matrix of the generators of~$A_q(G_2)$.
We  use
$\mathrm{Osc}_{q^3} =
\langle  {\rm{\bf A}}^+, {\rm{\bf A}}^-, {\rm{\bf K}} \rangle$
in addition to
$\mathrm{Osc}_q =
\langle  {\rm{\bf a}}^+, {\rm{\bf a}}^-, {\rm{\bf k}} \rangle$
(\ref{a:qoa}).
The fundamental representations
$\pi_i: A_q(G_2) \rightarrow \mathrm{End}(\mathcal{F}_{q_i})$
$(i=1,2)$ are given by{\samepage
\begin{gather}
 \pi_1(T)=
 \left(\!
\begin{array}{@{}c@{\,}c@{\,}c@{\,}c@{\,}c@{\,}c@{\,}c@{}}
\mu_1 {\rm{\bf a}}^- & \alpha_1 {\rm{\bf k}} &  0 & 0 & 0 & 0 &0 \\
-q\alpha^{-1}_1{\rm{\bf k}} & \mu^{-1}_1{\rm{\bf a}}^+ & 0 & 0&0 & 0& 0\\
0&0 & (\mu_1{\rm{\bf a}}^-)^2 &
[2]_1\alpha_1\mu_1 {\rm{\bf k}}\, {\rm{\bf a}}^- & (\alpha_1{\rm{\bf k}})^2 & 0& 0\\
0 & 0& -q\alpha_1^{-1}\mu_1{\rm{\bf a}}^- \,{\rm{\bf k}} &
{\rm{\bf a}}^- {\rm{\bf a}}^+ - {\rm{\bf k}}^2 & \alpha_1\mu^{-1}_1{\rm{\bf k}}
\,{\rm{\bf a}}^+ & 0& 0\\
0& 0& (q\alpha_1^{-1}{\rm{\bf k}})^2 &
-[2]_1(\alpha_1\mu_1)^{-1} {\rm{\bf k}}
\,{\rm{\bf a}}^+ & (\mu_1^{-1}{\rm{\bf a}}^+)^2 & 0& 0\\
0&0&0&0&0& \mu_1{\rm{\bf a}}^- & \alpha_1{\rm{\bf k}} \\
0&0&0&0&0& -q\alpha_1^{-1}{\rm{\bf k}} & \mu_1^{-1}{\rm{\bf a}}^+
 \end{array}
\!\right)
 , \nonumber\\
\pi_2(T)=
\begin{pmatrix}
1&  0&  0& 0& 0& 0&0\\
 0& \mu_2{\rm{\bf A}}^- & \alpha_2 {\rm{\bf K}} &0 &0 &0 &0 \\
 0& -q^3\alpha_2^{-1} {\rm{\bf K}} & \mu_2^{-1} {\rm{\bf A}}^+ & 0& 0& 0&0\\
0 &0 & 0& 1& 0&0 &0\\
 0& 0& 0& 0& \mu_2 {\rm{\bf A}}^- & \alpha_2 {\rm{\bf K}} & 0\\
0 & 0& 0& 0& -q^3\alpha^{-1}_2 {\rm{\bf K}} & \mu_2^{-1} {\rm{\bf A}}^+ & 0\\
 0&0&0&0&0&0& 1
 \end{pmatrix},\label{a:piTg}
\end{gather}
where $\alpha_i$, $\mu_i$ are nonzero parameters
and $[2]_1=q+q^{-1}$ as def\/ined after~(\ref{a:serre}).}

The Weyl group $W =\langle s_1, s_2 \rangle$ is the Coxeter system
with the relations
\begin{gather*}
s_1^2=s_2^2=1,\quad
s_2s_1s_2s_1s_2s_1 = s_1s_2s_1s_2s_1s_2.
\end{gather*}
Thus we have the isomorphism
$\pi_{212121} \simeq \pi_{121212}$.
Let $\Phi$ be the corresponding intertwiner and denote by $F$
the checked intertwiner $\Phi^\vee$
\begin{gather}
 \pi_{212121} \Phi = \Phi\pi_{121212},\qquad
\pi_{212121} F = F \pi'_{212121},\qquad
\pi'_{212121}
= P_{16}P_{25}P_{34}\pi_{121212}P_{16}P_{25}P_{34},\nonumber
\\
F = \Phi  P_{16}P_{25}P_{34}\in
\mathrm{End}(\mathcal{F}_{q^3}\otimes
\mathcal{F}_{q}\otimes
\mathcal{F}_{q^3}\otimes
\mathcal{F}_{q}\otimes
\mathcal{F}_{q^3}\otimes
\mathcal{F}_{q}).\label{a:fdef}
\end{gather}
Def\/ine the matrix elements of $F$ and its parameter-free part
$\mathscr{F}$ by
\begin{gather*}
F |i,j,k,l,m,n \rangle=
\sum_{a,b,c,d,e,f}F^{abcdef}_{ijklmn}|a,b,c,d,e,f \rangle,\\
F^{a b c d e f}_{i j k l m n} =
\mu_1^{3c-3k+d-l+3e-3m}
\mu_2^{2k-2c+l-d+3m-3e+n-f}
\mathscr{F}^{a b c d e f}_{i j k l m n}.
\end{gather*}
Then the following properties are valid for
$\mathscr{F} = (\mathscr{F}^{abcdef}_{ijklmn})$:
\begin{gather}
\mathscr{F}^{abcdef}_{ijklmn} \in \mathbb{Z}[q],\nonumber\\
\mathscr{F}^{abcdef}_{ijklmn}=0 \qquad \mathrm{unless}\ \
\left({a + b + 2c + d + e \atop b + 3c + 2d + 3e + f}\right)=
\left({i + j + 2k + l + m\atop j + 3k + 2l + 3m + n}\right),
\label{a:claw3}\\
\mathscr{F}^{-1} = \mathscr{F},\qquad
\mathscr{F}^{abcdef}_{ijklmn}
=\frac{(q^6)_i(q^2)_j(q^6)_k(q^2)_l(q^6)_m(q^2)_n}
{(q^6)_a(q^2)_b(q^6)_c(q^2)_d(q^6)_e(q^2)_f}
\mathscr{F}^{ijklmn}_{abcdef}.
\label{a:fprop}
\end{gather}
Due to (\ref{a:claw3}), $\mathscr{F}$ is the inf\/inite direct sum
of f\/inite-dimensional matrices.
The formula for $\mathscr{F}^{abcdef}_{ijklmn}|_{q=0}$
can be deduced by the ultradiscretization (tropical form) of \cite[Theorem~3.1(c)]{BZ97}.
Although a tedious algorithm can be formulated
for calculating any given $\mathscr{F}^{abcdef}_{ijklmn}$ by using~(\ref{a:fdef}),
an explicit formula for it is yet to be constructed.

\begin{example}
The following is the list of all the nonzero $\mathscr{F}^{abcdef}_{010101}$:
\begin{gather*}
\mathscr{F}^{000200}_{010101}
  = q^4 \big(1 - q^2\big) \big(1 - q^2 - q^4 - q^6\big),
\\
\mathscr{F}^{001001}_{010101}
  = -q \big(1 - q^2\big) \big(1 - q^2 - q^4 + q^8 + q^{10}\big), \\
\mathscr{F}^{010010}_{010101}
 = -q \big(1 - q^2\big) \big(1 - q^2 - q^4 + q^8 + q^{10}\big),
\\
\mathscr{F}^{010101}_{010101}
 =  1 - 2 q^2 + 2 q^6 + 3 q^8 - 2 q^{12} - 2 q^{14} - q^{16}, \\
\mathscr{F}^{020002}_{010101}
 =  q^4 \big({-}2 + 2 q^6 + q^8 + q^{10}\big),
\\
\mathscr{F}^{100011}_{010101}
 = -q^3 \big(1 - q^2\big) \big(1 - q^6 - q^8\big), \\
\mathscr{F}^{100102}_{010101}
 =  q \big(1 - q^2 - q^4 - q^6 + q^{10} + q^{12} + q^{14}\big),
\\
\mathscr{F}^{200004}_{010101}  = q^4, \\
\mathscr{F}^{11000 3}_{010101}
 = q \big(1 - q + q^2\big) \big(1 + q + q^2\big) \big(1 - q^2 - q^8\big).
\end{gather*}
\end{example}

\section{Main theorem}\label{a:sec5}

In this section we f\/ix two reduced words ${\bf i}=(i_1,\ldots,i_l)$,
${\bf j}=(j_1,\ldots,j_l)$ of the longest element $w_0\in W$.

\subsection[Definitions of $\gamma^A_B$ and $\Phi^A_B$]{Def\/initions of $\boldsymbol{\gamma^A_B}$ and $\boldsymbol{\Phi^A_B}$}

In the $U_q(\geh)$ side, we def\/ined
the PBW bases
$E_{\bf i}^A$, $E_{\bf j}^B$ of $U_q^+(\geh)$ in Section~\ref{a:ss:pbw}.
We def\/ine their transition coef\/f\/icient~$\gamma^A_B$ by
\begin{gather*}
E_{\bf i}^A=\sum_B \gamma^A_B E_{\bf j}^B.
\end{gather*}
While, in the $A_q(\geh)$ side, we have the intertwiner
$\Phi:
\mathcal{F}_{q_{i_1}}\otimes \cdots \otimes \mathcal{F}_{q_{i_l}}
\rightarrow
\mathcal{F}_{q_{j_1}}\otimes \cdots \otimes \mathcal{F}_{q_{j_l}}$
satisfying
\begin{gather}\label{a:piphi}
\pi_{\bf j}(g)\circ \Phi = \Phi \circ \pi_{\bf i}(g)\qquad
\forall\, g \in A_q(\mathfrak{g}).
\end{gather}
We take the parameters
$\mu$, $\alpha$ in (\ref{a:pisl2}) to be~1.
This in particular means for rank 2 cases that $\mu_i$, $\alpha_i$
entering $\pi_i(T)$ in (\ref{a:piTa}), (\ref{a:piTc}) and
(\ref{a:piTg}) are all~1.
The intertwiner $\Phi$ is normalized
by $\Phi |0,0,\ldots, 0\rangle =  |0,0,\ldots, 0\rangle$.
Under these conditions a matrix element~$\Phi^A_B$ of~$\Phi$
is uniquely specif\/ied by
\begin{gather*}
\Phi |B\rangle=\sum_A \Phi^A_B |A\rangle,
\end{gather*}
where $A=(a_1, \ldots, a_l) \in (\Z_{\ge 0})^l$ and
$|A\rangle = |a_1\rangle \otimes \cdots \otimes |a_l\rangle
\in \mathcal{F}_{q_{j_1}}\otimes \cdots \otimes \mathcal{F}_{q_{j_l}}$
and similarly for $|B\rangle
\in \mathcal{F}_{q_{i_1}}\otimes \cdots \otimes \mathcal{F}_{q_{i_l}}$.
Then our main result is
\begin{theorem}\label{theorem:gamma=R}
\[
\gamma^A_B=\Phi^A_B.
\]
\end{theorem}

For any pair $({\bf i}$, ${\bf j})$, from ${\bf i}$ one can reach~${\bf j}$ by applying Coxeter relations.
In view of the uniqueness of $\gamma$ and $\Phi$
and the fact that
the braid group action $T_i$ is an algebra homomorphism,
the proof of this theorem reduces to
establishing the same equality for all $\geh$ of rank~2.
This will  be done in the rest of this section.

\subsection{Proof of Theorem \ref{theorem:gamma=R} for rank 2 cases}

In the rank 2 cases,
there are two reduced expressions
$s_{i_1}\cdots s_{i_l}$ for the longest element of the
Weyl group.
Denote the associated sequences
${\bf i} = (i_1,\ldots, i_l)$  by ${\bf 1}$, ${\bf 2}$  and
set ${\bf 1}' = {\bf 2}$, ${\bf 2}' = {\bf 1}$.
Concretely, we take them as
\begin{alignat*}{5}
& A_2:  \quad && {\bf 1}=(1,2,1),\qquad &&{\bf 2}=(2,1,2), \qquad &&(q_1,q_2)=(q,q),& \\
& C_2: \quad && {\bf 1}=(1,2,1,2),\qquad &&{\bf 2}=(2,1,2,1), \qquad &&(q_1,q_2)=\big(q,q^2\big),& \\
& G_2: \quad && {\bf 1}=(1,2,1,2,1,2),\qquad &&{\bf 2}=(2,1,2,1,2,1), \qquad &&(q_1,q_2)=\big(q,q^3\big),&
\end{alignat*}
where $q_i$ def\/ined after (\ref{a:serre}) is also recalled.
In order to simplify the formulae in
Section \ref{a:ss:rank2}, we use the PBW bases and the
Fock states in yet another normalization as follows:
\begin{gather}
{\tilde E}^A_{\bf i}  := ([a_1]_{i_1}!\cdots [a_l]_{i_l}!) E^A_{\bf i}
= e^{a_1}_{\beta_1}\cdots e^{a_l}_{\beta_l}, \nonumber\\
|A\rangle\!\rangle  := d_{i_1,a_1}\cdots d_{i_l,a_l}|A\rangle,\qquad
d_{i,a} = q^{-a(a-1)/2}_i\lambda^a_i,\qquad
\lambda_i = \big(1-q_i^2\big)^{-1},\label{a:kett}
\end{gather}
where $A=(a_1,\ldots, a_l)$.
See after (\ref{a:serre}) for the symbol $[a]_i!$.
$e_{\beta_r}$ is def\/ined in (\ref{a:ebdef}).
Accordingly we introduce the matrix elements ${\tilde \gamma}^A_B$ and
${\tilde \Phi}^A_B$ by
\begin{gather*}
{\tilde E}^A_{\bf i} =\sum_B {\tilde \gamma}^A_B
{\tilde E}^B_{{\bf i}'},\qquad
\Phi |B\rangle\!\rangle=\sum_A {\tilde \Phi}^A_B |A\rangle\!\rangle,
\qquad {\bf i}={\bf 1}, {\bf 2}.
\end{gather*}
It follows that
$\gamma^A_B =
{\tilde \gamma}^A_B\prod\limits_{k=1}^l([b_k]_{i_k}!/[a_k]_{i_k}!)$
and
$\Phi^A_B=
{\tilde \Phi}^A_B\prod\limits_{k=1}^l(d_{i_k,a_k}/d_{i_k,b_k})$
for $B=(b_1,\ldots, b_l)$.
On the other hand, we know
$\Phi^A_B=\Phi^B_A\prod\limits_{k=1}^l((q_{i_k}^2)_{b_k}/(q_{i_k}^2)_{a_k})$
from (\ref{a:rprop}), (\ref{a:kprop}) and (\ref{a:fprop}).
Due to the identity
$(q_i^2)_md_{i,m}=[m]_i!$,
the assertion $\gamma^A_B=\Phi^A_B$
of Theorem \ref{theorem:gamma=R} is equivalent to
\begin{gather}\label{a:gyaku}
{\tilde \gamma}^A_B={\tilde \Phi^B_A}.
\end{gather}

Let $\rho_{\bf i}(x)=(\rho_{\bf i}(x)_{A B})$ be the matrix
for the left multiplication of  $x \in U^+_q(\geh)$:
\begin{gather}\label{a:lmul}
x \cdot {\tilde E}^{A}_{\bf i} = \sum_B  {\tilde E}^B_{\bf i}
\rho_{\bf i}(x)_{BA}.
\end{gather}
Let further $\pi_{\bf i}(g)=(\pi_{\bf i}(g)_{A B})$ be the
representation matrix of $g \in A_q(\geh)$:
\begin{gather}\label{a:pib}
\pi_{\bf i}(g) | A\rangle\!\rangle
= \sum_B |B\rangle\!\rangle \pi_{\bf i}(g)_{BA}.
\end{gather}

The following element in
the right quotient ring $A_q(\geh)_\mathcal{S}$
will play a key role in our proof.
\begin{gather}\label{a:xidef}
\xi_i = \lambda_i (\sigma_i e_i)/\sigma_i,\qquad i=1,2.
\end{gather}
See Def\/inition \ref{def:sigma} for $\sigma_i$
and (\ref{a:kobetsu}),  (\ref{c:kobetsu}),  (\ref{g:kobetsu})
for the concrete forms in rank 2 cases.
In Section \ref{a:ss:rank2} we will check
the following statement case by case.
\begin{proposition}\label{a:pr:rnk2}
For $\geh$ of rank $2$,
$\pi_{\bf i}(\sigma_i)$ is invertible and
the following equality is valid:
\begin{gather}\label{a:key}
\rho_{\bf i}(e_i)_{A B} = \pi_{\bf i}(\xi_i)_{A B}, \qquad i=1,2,
\end{gather}
where the right-hand side means
$\lambda_i\pi_{\bf i}(\sigma_ie_i) \pi_{\bf i}(\sigma_i)^{-1}$.
\end{proposition}

\begin{proof}[Proof of Theorem \ref{theorem:gamma=R} for rank 2 case.]
We write the both sides of (\ref{a:key}) as $M^i_{AB}$ and the one
for ${\bf i}'$ instead of ${\bf i}$ as $M^{\prime i}_{AB}$.
From
\begin{gather*}
\sum_{B,C} {\tilde E}^C_{{\bf i}'} M^{\prime i}_{CB}{\tilde \gamma}^A_B=
e_i \sum_B {\tilde E}^B_{{\bf i}'}{\tilde \gamma}^A_B =
e_i {\tilde E}^A_{\bf i} = \sum_B {\tilde E}^{B}_{\bf i}  M^i_{BA}
= \sum_{B,C}{\tilde E}^C_{{\bf i}'}
{\tilde \gamma}^B_C M^i_{BA}
\end{gather*}
we have
$\sum_{B} M^{\prime i}_{CB}{\tilde \gamma}^A_B
= \sum_{B}{\tilde \gamma}^B_CM^i_{BA}$.
On the other hand, the action of the two sides of (\ref{a:piphi}) with
$g=\xi_i$ and
${\bf j}={\bf i}'$ are calculated as
\[
\pi_{{\bf i}'}(\xi_i)\circ \Phi|A\rangle\!\rangle
=\pi_{{\bf i}'}(\xi_i)\sum_B|B\rangle\!\rangle{\tilde \Phi}^B_A
=\sum_{B,C}|C\rangle\!\rangle
M_{CB}^{\prime i}{\tilde \Phi}_A^B
\]
and
\[
\Phi \circ \pi_{{\bf i}}(\xi_i)|A\rangle\!\rangle
=\Phi\sum_B|B\rangle\!\rangle M_{BA}^i
=\sum_{B,C}|C\rangle\!\rangle {\tilde \Phi}_B^C M_{BA}^i.
\]
Hence
$\sum_{B}
M_{CB}^{\prime i}{\tilde \Phi}_A^B
= \sum_{B}{\tilde \Phi}_B^C M_{BA}^i$.
Thus ${\tilde \gamma}^A_B$ and ${\tilde \Phi}^B_A$
satisfy the same relation.
Moreover the maps $\pi_{\bf i}$ and
$\rho_{\bf i}$ are both homomorphism, i.e.,
$\pi_{\bf i}(gh)=\pi_{\bf i}(g)\pi_{\bf i}(h)$
and
$\rho_{\bf i}(xy)=\rho_{\bf i}(x)\rho_{\bf i}(y)$.
We know that ${\Phi}$ is the intertwiner of the
irreducible $A_q(\geh)$ modules and
(\ref{a:gyaku}) obviously holds as $1=1$ at
$A=B=(0,\ldots,0)$.
Thus it is valid for arbitrary $A$ and $B$.
\end{proof}

\begin{conjecture}
The equality \eqref{a:key} is valid for any $\geh$.
\end{conjecture}

\subsection{Explicit formulae for rank 2 cases:
Proof of Proposition~\ref{a:pr:rnk2}}\label{a:ss:rank2}

Here we present the explicit formulae of
(\ref{a:lmul}) with $x=e_i$
and (\ref{a:pib}) with $g=\sigma_i, \sigma_ie_i$ that allow one to check
Proposition \ref{a:pr:rnk2}.
We use the notation $\qmq{i}=q^i-q^{-i}$.
In each case, there are two ${\bf i}$-sequences,
${\bf 1}$ and ${\bf 2}={\bf 1}'$
corresponding to the two reduced words.
Let $\chi$ be the anti-algebra involution
such that $\chi(e_i)=e_i$.
Then the relation
$\chi({\tilde E}^{A}_{\bf i}) = {\tilde E}^{\bar A}_{{\bf i}'}$ holds,
where ${\bar A} = (a_l,\ldots, a_2,a_1)$ denotes the reversal of
$A=(a_1,a_2,\ldots, a_l)$.
Applying $\chi$ to (\ref{a:lmul}) with $x=e_i$ yields the
right multiplication formula
${\tilde E}^{\bar A}_{{\bf i}'}\cdot e_i  = \sum_B  {\tilde E}^{\bar B}_{{\bf i}'}
\rho_{\bf i}(e_i)_{BA}$ for ${\bf i}'$-sequence.
In view of this fact, we shall present
the left and right multiplication formulae
for ${\bf i}={\bf 2}$ only.

As for (\ref{a:pib})  with $g=\xi_i$ in (\ref{a:xidef}),
explicit formulae for $\sigma_i, \sigma_i e_i \in A_q(\geh)$ and
their image by the both representations
$\pi_{\bf 1}$ and $\pi_{\bf 2}$ will be given.
We include an exposition on how to use these data to check (\ref{a:key})
along the simplest $A_2$ case.
The $C_2$ and $G_2$ cases are similar.

\subsubsection[$A_2$ case]{$\boldsymbol{A_2}$ case}


\centerline{\includegraphics{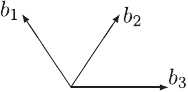}}

The $q$-Serre relations are
\begin{gather*}
e_1^2 e_2-[2]_1e_1 e_2 e_1+e_2 e_1^2=0,\qquad
e_2^2 e_1-[2]_1e_2 e_1 e_2+e_1 e_2^2=0,
\end{gather*}
where $[m]_1=\langle m\rangle/\langle 1\rangle$.
Let $b_1$, $b_2$, $b_3$ be the generator for positive roots: $b_1=e_2$,
$b_2=e_1 e_2-q e_2 e_1$ and $b_3=e_1$.
In the notation of Section~\ref{a:ss:pbw},
they are the root vectors
$b_i = e_{\beta_i} $ associated with the reduced expression
$w_0 = s_2s_1s_2$ for ${\bf 2}=(2,1,2)$.
The corresponding
positive roots are
$(\beta_1, \beta_2, \beta_3)=(\alpha_2, \alpha_1+\alpha_2, \alpha_1)$.
In particular, $b_2 = T_2(e_1)$.
Their commutation relations are
$b_2  b_1 = q^{-1}{b_1 b_2}$,
$b_3  b_1 = b_2+q b_1 b_3$,
$b_3  b_2 = q^{-1}{b_2 b_3}$.

\begin{lemma} \label{lemma:rho-rep-A}
For $\undividedB{a,b,c}=b_1^a b_2^b b_3^c$, we have
\begin{gather*}
\undividedB{a,b,c} \cdot e_1 =\undividedB{a,b,c+1},\\
\undividedB{a,b,c} \cdot e_2 =q^{c-b} \undividedB{a+1,b,c}+[c]_1 \undividedB{a,b+1,c-1},\\
e_1 \cdot \undividedB{a,b,c} =q^{a-b} \undividedB{a,b,c+1}+[a]_1 \undividedB{a-1,b+1,c},\\
e_2 \cdot \undividedB{a,b,c} =\undividedB{a+1,b,c}.
\end{gather*}
\end{lemma}

\begin{proof}
By induction, we have
\begin{gather*}
{b_3 b_1^n} = q^n {b_1^n b_3}+[n]_1 {b_1^{n-1} b_2},  \qquad
{b_3 b_2^n} = q^{-n} {b_2^n b_3},\\
{b_3^n b_1} = q^n {b_1 b_3^n}+[n]_1 {b_2 b_3^{n-1}},  \qquad
{b_2^n b_1} = q^{-n} {b_1 b_2^n}.
\end{gather*}
The lemma is a direct consequence of these formulae.
\end{proof}

Set
${\tilde E}^{a,b,c}_{\bf 1} = \chi({\tilde E}^{c,b,a}_{\bf 2})
= \chi(b^a_3)\chi(b^b_2)\chi(b^c_1)
=b_3^ab_2^{\prime b}b_1^c$,
where $b'_2 := \chi(b_2) = e_2e_1-qe_1e_2$.
By applying $\chi$ to the f\/irst two relations
in Lemma~\ref{lemma:rho-rep-A}, we get
\begin{gather}\label{a:kore}
e_1\cdot {\tilde E}^{a,b,c}_{\bf 1} = {\tilde E}^{a+1,b,c}_{\bf 1},
\qquad
e_2\cdot {\tilde E}^{a,b,c}_{\bf 1} = q^{a-b}E^{a,b,c+1}_{\bf 1}
+[a]_1{\tilde E}^{a-1,b+1,c}_{\bf 1}.
\end{gather}
Thus we f\/ind $\rho_{{\bf i}'}(e_i) = \rho_{\bf i}(e_{3-i})$.
This property is only valid for $A_2$ and not in $C_2$ and $G_2$.

Let us turn to the representations $\pi_{{\bf i}}$ of $A_q(A_2)$.
The elements $\sigma_i$ in Def\/inition \ref{def:sigma}
and $\sigma_ie_i$
are given by{\samepage
\begin{gather}\label{a:kobetsu}
\sigma_1 = t_{13}, \qquad
\sigma_2 = t_{12}t_{23}-qt_{22}t_{13},\qquad
\sigma_1 e_1 = t_{23},\qquad
\sigma_2e_2 = t_{12}t_{33}-qt_{32}t_{13}.
\end{gather}
See the exposition at the end of Section \ref{a:sss:a1} and the remark after
Lemma \ref{le:ff}.}

From (\ref{a:cpro}) and (\ref{a:piTa}) with $\alpha_i=\mu_i=1$, we f\/ind
\begin{gather*}
\pi_{{\bf 1}}(\sigma_1) = {\bf k}_1{\bf k}_2,\qquad\!
\pi_{{\bf 1}}(\sigma_1e_1)
= {\rm{\bf a}}^+_1{\bf k}_2,\qquad\!
\pi_{{\bf 1}}(\sigma_2) ={\bf k}_2{\bf k}_3,\qquad\!
\pi_{{\bf 1}}(\sigma_2e_2) =
{\rm{\bf a}}^-_1{\rm{\bf a}}^+_2{\bf k}_3
+{\bf k}_1{\rm{\bf a}}^+_3,
\end{gather*}
where the notation like
${\bf k}_1{\rm{\bf a}}^+_3 = {\bf k}\otimes 1 \otimes {\rm{\bf a}}^+$
has been used.
Since ${\bf k}\in \mathrm{End}({\mathcal F}_q)$ is invertible,
so is $\pi_{\bf i}(\sigma_i)$ and we may write
\begin{gather*}
\pi_{\bf 1}(\xi_1) = \lambda_1
{\rm{\bf a}}^+_1{\bf k}^{-1}_1,\qquad
\pi_{\bf 1}(\xi_2) = \lambda_2
({\rm{\bf a}}^-_1{\rm{\bf a}}^+_2{\bf k}^{-1}_2
+{\bf k}_1{\bf k}^{-1}_2{\rm{\bf a}}^+_3{\bf k}^{-1}_3),
\end{gather*}
where $\lambda_1=\lambda_2 = (1-q^2)^{-1}$.
The action of each component on the ket vector
$|m\rangle\!\rangle :=d_{i,m}|m\rangle \in {\mathcal F}_{q_i}$
(cf.~(\ref{a:kett})) takes the form
\begin{gather}\label{a:ketta}
{\rm{\bf a}}^+|m\rangle\!\rangle =
\lambda_i^{-1}q_i^m|m+1\rangle\!\rangle,\qquad
{\rm{\bf a}}^-|m\rangle\!\rangle =
[m]_i|m-1\rangle\!\rangle,\qquad
{\bf k}|m\rangle\!\rangle=q_i^m|m\rangle\!\rangle,
\end{gather}
due to (\ref{a:akf}).
(The formula \eqref{a:ketta} is valid also for $C_2$ and $G_2$
provided that ${\rm{\bf a}}^+$, ${\rm{\bf a}}^-$, ${\rm{\bf k}}$
are interpreted as
${\rm{\bf A}}^+$, ${\rm{\bf A}}^-$, ${\rm{\bf K}}$
for $i=2$.)
Thus one has
\begin{gather*}
\pi_{\bf 1}(\xi_1)|a,b,c\rangle\!\rangle =
|a+1,b,c\rangle\!\rangle,
\qquad
\pi_{\bf 1}(\xi_2)|a,b,c\rangle\!\rangle =
[a]_1|a-1,b+1,c\rangle\!\rangle
+q^{a-b}|a,b,c+1\rangle\!\rangle.
\end{gather*}
This agrees with (\ref{a:kore}) thereby
proving (\ref{a:key}) for ${\bf i}={\bf 1}$.
The other case  ${\bf i}={\bf 2}$ also holds due to the
symmetry $\pi_{\bf 2}(\xi_i) = \pi_{\bf 1}(\xi_{3-i})$.
Thus Proposition \ref{a:pr:rnk2} is established for~$A_2$.

In terms of the checked intertwiner $\mathscr{R}$ in Section~\ref{a:ss:a2},
Theorem~\ref{theorem:gamma=R} implies
\begin{gather*}
E^{a,b,c}_{\bf i} = \sum_{i,j,k}
\mathscr{R}^{a b c}_{i j k} E^{k,j,i}_{{\bf i}'}.
\end{gather*}
This is valid either for ${\bf i}={\bf 1}$ or ${\bf 2}$
thanks to the middle property in~(\ref{a:rprop}).
This relation connecting the PBW bases with
the solution of the tetrahedron equation is due to~\cite{S08}.

\subsubsection[$C_2$ case]{$\boldsymbol{C_2}$ case}


\centerline{\includegraphics{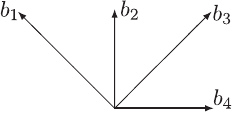}}

The $q$-Serre relations are
\begin{gather*}
 \eshort^3 \elong
-[3]_1\eshort^2 \elong \eshort
+[3]_1\eshort \elong \eshort^2
-\elong \eshort^3=0,\qquad
 \elong^2 \eshort
-[2]_2\elong \eshort \elong
+\eshort \elong^2=0,
\end{gather*}
where $[m]_1 = \qmq{m}/\qmq{1}$
and $[m]_2 = \qmq{2m}/\qmq{2}$.

Let $b_1, \dots, b_4$ be the generator for positive  roots:
$b_1=\elong$,
$b_2=\eshort\elong-q^{2} \elong\eshort$,
$b_3=\frac{1}{[2]_1}(\eshort b_2-b_2 \eshort)$
and $b_4=\eshort$.
Their commutation relations are
$b_2  b_1 = q^{-2} b_1 b_2$,
$b_3  b_1 = -q^{-1} \qmq{1}[2]_1^{-1} b_2^2+b_1 b_3$,
$b_4  b_1 = b_2+q^{2} b_1 b_4$,
$b_3  b_2 = q^{-2} b_2 b_3$,
$b_4  b_2 = [2]_1 b_3+b_2 b_4$,
$b_4  b_3 = q^{-2} b_3 b_4$.

\begin{lemma} \label{lemma:rho-rep-C}
For $\undividedB{a,b,c,d}=b_1^a b_2^b b_3^c b_4^d$, we have
\begin{gather*}
\undividedB{a,b,c,d} \cdot \eshort
 =\undividedB{a,b,c,d+1},\\
\undividedB{a,b,c,d} \cdot \elong
 =[d]_1 q^{d-2 c-1} \undividedB{a,b+1,c,d-1}+q^{2 (d-b)} \undividedB{a+1,b,c,d}\\
\hphantom{\undividedB{a,b,c,d} \cdot \elong=}{}
-\qmq{1} q^{2 d-2 c+1}[c]_2 [2]_1^{-1} \undividedB{a,b+2,c-1,d}+[d-1]_1 [d]_1
\undividedB{a,b,c+1,d-2},\\
\eshort \cdot \undividedB{a,b,c,d}
 = [2]_1 [b]_1 q^{2 a-b+1} \undividedB{a,b-1,c+1,d}
+q^{2 a-2 c} \undividedB{a,b,c,d+1}+[a]_2 \undividedB{a-1,b+1,c,d},\\
\elong  \cdot \undividedB{a,b,c,d}
=\undividedB{a+1,b,c,d}.
\end{gather*}
\end{lemma}

\begin{proof}
 By induction, we have
\begin{gather*}
 {b_4 b_1^n} = {b_1^n b_4} q^{2 n}+[n]_2 {b_1^{n-1}b_2},\qquad
 {b_4 b_2^n} = [2]_1 [n]_1 {b_2^{n-1} b_3} q^{-n+1}+{b_2^n b_4},\\
 {b_4 b_3^n} = q^{-2 n} {b_3^n b_4},\qquad
 {b_4^n b_1} = [n]_1 {b_2 b_4^{n-1}} q^{n-1}+{b_1 b_4^n} q^{2 n}+[n-1]_1 [n]_1 {b_3 b_4^{n-2}},\\
 {b_3^n b_1} = -q^{1-2 n}\qmq{1} [n]_2 [2]_1^{-1} {b_2^2 b_3^{n-1}}+{b_1 b_3^n},\qquad
 {b_3^n b_2} = q^{-2 n} {b_2 b_3^n},\qquad
 {b_2^n b_1} = q^{-2 n} {b_1 b_2^n}.
 \end{gather*}
The lemma is a direct consequence of these formulae.
\end{proof}

Set ${\tilde E}_{\bf 1}^{a,b,c,d} = \chi\big({\tilde E}^{d,c,b,a}_{\bf 2}\big)$.
The left multiplication formula for this basis is deduced from
the above lemma by applying~$\chi$.
One can adjust the def\/inition of
$E^A_{\bf i}$ in~(\ref{a:Eai}) with that in~\cite{Xi} by setting $v=q^{-1}$.

Let us turn to the representations $\pi_{{\bf i}}$ of $A_q(C_2)$.
The elements $\sigma_i$ in Def\/inition~\ref{def:sigma}
and~$\sigma_ie_i$
are given by
\begin{gather}\label{c:kobetsu}
\sigma_1 = t_{14}, \qquad
\sigma_2 = t_{13}t_{24}-qt_{23}t_{14},\qquad
\sigma_1 e_1 = t_{24},\qquad
\sigma_2e_2 = t_{13}t_{34}-qt_{33}t_{14}.
\end{gather}
From (\ref{a:cpro}) and (\ref{a:piTc}) with $\alpha_i=\mu_i=1$, we have
\begin{gather*}
\pi_{\bf 1}(\sigma_1)=- {{\bf k}_1}
{{\bf K}_2} {{\bf k}_3},\\
\pi_{\bf 1}(\sigma_1e_1)=- {{\bf a}^+_1}
{{\bf K}_2} {{\bf k}_3},\\
\pi_{\bf 1}(\sigma_2)=- {{\bf K}_2}
{{\bf k}_3}^{2} {{\bf K}_4},\\
\pi_{\bf 1}(\sigma_2e_2)=- {{\bf a}^-_1}^2
{{\bf A}^+_2} {{\bf k}_3}^{2} {{\bf K}_4}
-[2]_1 {{\bf a}^-_1} {{\bf k}_1}
{{\bf a}^+_3} {{\bf k}_3} {{\bf K}_4}
- {{\bf k}_1}^{2} {{\bf A}^-_2} {{\bf a}^+_3}^2
{{\bf K}_4}
- {{\bf A}^+_4} {{\bf k}_1}^{2} {{\bf K}_2},\\
\lambda_1^{-1}\pi_{\bf 1}(\xi_1)={{\bf a}^
+_1} {{\bf k}_1}^{-1},\\
\lambda_2^{-1}\pi_{\bf 1}(\xi_2) ={{\bf a}^-
_1}^2 {{\bf A}^+_2} {{\bf K}_2}^{-1}
+{{\bf k}_1}^{2} {{\bf A}^-_2} {{\bf K}_2}^{-1}
{{\bf a}^+_3}^2 {{\bf k}_3}^{-2}\\
\hphantom{\lambda_2^{-1}\pi_{\bf 1}(\xi_2) =}{}
+[2]_1 {{\bf a}^-_1} {{\bf k}_1}
{{\bf K}_2}^{-1} {{\bf a}^+_3} {{\bf k}_3}^{-1}
+{{\bf k}_1}^{2} {{\bf k}_3}^{-2} {{\bf A}^+_4}
{{\bf K}_4}^{-1},\\
\pi_{\bf 2}(\sigma_1)=- {\bf k}_2 {\bf
K}_3 {\bf k}_4,\\
\pi_{\bf 2}(\sigma_1e_1) =- {{\bf K}_1}
{\bf k}_2 {{\bf a}^+_4}
- {{\bf K}_1} {{\bf a}^-_2} {{\bf A}^+_3} {\bf
k}_4
- {{\bf A}^-_1} {{\bf a}^+_2} {\bf K}_3 {\bf k}
_4,\\
\pi_{\bf 2}(\sigma_2)=- {{\bf K}_1}
{\bf k}_2^2 {\bf K}_3,\\
\pi_{\bf 2}(\sigma_2e_2)=- {{\bf A}^+_1}
{\bf k}_2^2 {\bf K}_3,\\
\lambda_1^{-1}\pi_{\bf 2}(\xi_1)={{\bf A}^-
_1} {{\bf a}^+_2} {\bf k}_2^{-1}
+{{\bf K}_1} {{\bf a}^-_2} {\bf k}_2^{-1}
{{\bf A}^+_3} {\bf K}_3^{-1}
+{{\bf K}_1} {\bf K}_3^{-1} {{\bf a}^+_4} {\bf
k}_4^{-1},\\
\lambda_2^{-1}\pi_{\bf 2}(\xi_2)={{\bf A}^
+_1} {{\bf K}_1}^{-1}.
\end{gather*}

We f\/ind that $\pi_{\bf i}(\sigma_i)$ is invertible.
Comparing these formulae with Lemma~\ref{lemma:rho-rep-C}
by using (\ref{a:ketta}), the equality (\ref{a:key}) is directly checked.
Thus Proposition \ref{a:pr:rnk2} is established for~$C_2$.

In terms of the checked intertwiner $\mathscr{K}$ in Section~\ref{a:ss:c2},
Theorem~\ref{theorem:gamma=R} implies
\begin{gather*}
E^{a,b,c,d}_{\bf 2} = \sum_{i,j,k,l}
\mathscr{K}^{a b c d}_{i j k l} E^{l,k,j,i}_{\bf 1}.
\end{gather*}
Thus the solution to the 3D ref\/lection equation~\cite{KO1} is
identif\/ied with the transition coef\/f\/icient of the PBW bases for
$U^+_q(C_2)$.

\subsection[$G_2$ case]{$\boldsymbol{G_2}$ case}


\centerline{\includegraphics{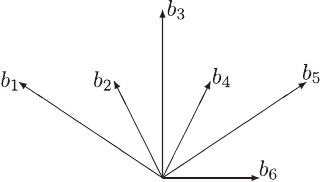}}

The $q$-Serre relations are
\begin{gather*}
 \eshort^4 \elong
-[4]_1\eshort^3 \elong \eshort
+[4]_1[3]_1/[2]_1^{-1}\eshort^2 \elong \eshort^2
-[4]_1\eshort \elong \eshort^3
+\elong \eshort^4=0,\\
 \elong^2 \eshort
-[2]_2\elong \eshort \elong
+\eshort \elong^2=0,
\end{gather*}
where we remind that $[m]_1 = \qmq{m}/\qmq{1}$
and $[m]_2 = \qmq{3m}/\qmq{3}$.

Let $b_1, \dots, b_6$ be the generator for positive  roots:
$b_1=\elong$,
$b_2=\eshort \elong-q^3 \elong \eshort$,
$b_4=\frac{1}{[2]_1}(\eshort b_2-q b_2 \eshort)$,
$b_5=\frac{1}{[3]_1}(\eshort b_4-q^{-1} b_4 \eshort)$,
$b_3=\frac{1}{[3]_1}(b_4 b_2-q^{-1} b_2 b_4)$
and $b_6=\eshort$.
Their commutation relations are as follows:
$b_2 b_1 = b_1 b_2 q^{-3}$,
$b_3 b_1 = \qmq{1}^2 b_2^3 q^{-3} [3]_1^{-1}+b_1 b_3 q^{-3}$,
$b_4 b_1 = b_1 b_4-b_2^2 \qmq{1} q^{-1}$,
$b_5 b_1 = b_1 b_5 q^3-b_2 b_4 \qmq{1} q^{-1}-(q^4+q^2-1) b_3 q^{-3}$,
$b_6 b_1 = b_1 b_6 q^3+b_2$,
$b_3 b_2 = b_2 b_3 q^{-3}$,
$b_4 b_2 = b_2 b_4 q^{-1}+b_3 [3]_1$,
$b_5 b_2 = b_2 b_5-b_4^2 \qmq{1} q^{-1}$,
$b_6 b_2 = q b_2 b_6+b_4 [2]_1$,
$b_4 b_3 = b_3 b_4 q^{-3}$,
$b_5 b_3 = \qmq{1}^2 b_4^3 q^{-3} [3]_1^{-1}+b_3 b_5 q^{-3}$,
$b_6 b_3 = b_3 b_6-b_4^2 \qmq{1} q^{-1}$,
$b_5 b_4 = b_4 b_5 q^{-3}$,
$b_6 b_4 = [3]_1 b_5+b_4 b_6 q^{-1}$,
$b_6 b_5 = b_5 b_6 q^{-3}$.

\begin{lemma} \label{lemma:rho-rep-G}
For $\undividedB{a,b,c,d,e,f}=b_1^a b_2^b\cdots b_6^f$, we have
\begin{gather*}
\undividedB{a,b,c,d,e,f} \cdot \eshort=
\undividedB{a,b,c,d,e,f+1},\\
\undividedB{a,b,c,d,e,f} \cdot \elong=
 -\qmq{1} [e]_2 q^{-3 c-d+3 f-1} \undividedB{a,b+1,c,d+1,e-1,f} \\
 \hphantom{\undividedB{a,b,c,d,e,f} \cdot \elong=}{}
 +\qmq{1}^2 [e-1]_2 [e]_2 [3]_1^{-1}q^{-3 e+3 f+3} \undividedB{a,b,c,d+3,e-2,f}\\
 \hphantom{\undividedB{a,b,c,d,e,f} \cdot \elong=}{}
 -\qmq{3} [d-1]_1 [d]_1 q^{-3 c-2 d+3 e+3 f+1} \undividedB{a,b+1,c+1,d-2,e,f}\\
 \hphantom{\undividedB{a,b,c,d,e,f} \cdot \elong=}{}
 -\qmq{1} [d]_1 q^{-6 c-d+3 (e+f)} \undividedB{a,b+2,c,d-1,e,f} \!
  +[f\!-\!1]_1 [f]_1 q^{-3 e+f-2} \undividedB{a,b,c,d+1,e,f-2}\\
\hphantom{\undividedB{a,b,c,d,e,f} \cdot \elong=}{}
 +[3]_1 [d]_1 [f]_1 q^{2 f-2 d} \undividedB{a,b,c+1,d-1,e,f-1}
 +[f]_1 q^{-3 c-d+2 f-2} \undividedB{a,b+1,c,d,e,f-1}\\
\hphantom{\undividedB{a,b,c,d,e,f} \cdot \elong=}{}
 +q^{-3 (b+c-e-f)} \undividedB{a+1,b,c,d,e,f}
 +\qmq{1}^2 [c]_2 [3]_1^{-1}q^{3 (-2 c+e+f+1)} \undividedB{a,b+3,c-1,d,e,f}\\
\hphantom{\undividedB{a,b,c,d,e,f} \cdot \elong=}{}
 -\qmq{3} [d-2]_1 [d-1]_1 [d]_1 q^{3 (-d+e+f+2)} \undividedB{a,b,c+2,d-3,e,f} \\
\hphantom{\undividedB{a,b,c,d,e,f} \cdot \elong=}{}
 -\qmq{1} [e]_2 [f]_1 q^{-3 e+2 f} \undividedB{a,b,c,d+2,e-1,f-1}\\
\hphantom{\undividedB{a,b,c,d,e,f} \cdot \elong=}{}
 -[e]_2 q^{-3 d+3 f} (q^{2d+1} [3]_1-[2]_2) \undividedB{a,b,c+1,d,e-1,f} \\
\hphantom{\undividedB{a,b,c,d,e,f} \cdot \elong=}{}
 +[f-2]_1 [f-1]_1 [f]_1 \undividedB{a,b,c,d,e+1,f-3},\\
\eshort \cdot \undividedB{a,b,c,d,e,f}=
 -\qmq{1} [c]_2 q^{3 a+b-3 c+2} \undividedB{a,b,c-1,d+2,e,f}
 +[3]_1 [b-1]_1 [b]_1 q^{3 a-b+2} \undividedB{a,b-2,c+1,d,e,f}\\
\hphantom{\eshort \cdot \undividedB{a,b,c,d,e,f}=}{}
 +[3]_1 [d]_1 q^{3 a+b-2 d+2} \undividedB{a,b,c,d-1,e+1,f}
 +q^{3 a+b-d-3 e} \undividedB{a,b,c,d,e,f+1}\\
\hphantom{\eshort \cdot \undividedB{a,b,c,d,e,f}=}{}
 +[2]_1 [b]_1 q^{3 (a-c)} \undividedB{a,b-1,c,d+1,e,f}
 +[a]_2 \undividedB{a-1,b+1,c,d,e,f},\\
 \elong \cdot \undividedB{a,b,c,d,e,f}=
\undividedB{a+1,b,c,d,e,f}.
\end{gather*}
\end{lemma}
\begin{proof} By induction, we have
\begin{gather*}
 b_6 b_1^n = q^{3 n} b_1^n b_6+[n]_2 b_1^{n-1} b_2,\\
 b_6 b_2^n = [3]_1 q^{2-n} [n-1]_1 [n]_1 b_2^{n-2} b_3+q^n b_2^n b_6+[2]_1 [n]_1 b_2^{n-1} b_4,\\
 b_4 b_3^n = q^{-3 n} b_3^n b_4,\qquad
 b_6 b_3^n = b_3^n b_6-\qmq{1} q^{2-3 n} [n]_2 b_3^{n-1} b_4^2,\\
 b_6 b_4^n = [3]_1 q^{2-2 n} [n]_1 b_4^{n-1} b_5+q^{-n} b_4^n b_6, \qquad
 b_6 b_5^n = q^{-3 n} b_5^n b_6,
 \end{gather*}
and
\begin{gather*}
 b_6^n b_1 = q^{n-2} [n-1]_1 [n]_1 b_4 b_6^{n-2}+q^{3 n} b_1 b_6^n \\
\hphantom{b_6^n b_1 =}{}
 +q^{2 (n-1)} [n]_1 b_2 b_6^{n-1}+[n-2]_1 [n-1]_1 [n]_1 b_5 b_6^{n-3},\\
 b_5^n b_1 = \qmq{1}^2 q^{-3 (n-1)} [n-1]_2 [n]_2 [3]_1^{-1}b_4^3 b_5^{n-2}+q^{3 n} b_1 b_5^n\\
\hphantom{b_5^n b_1 =}{}
 -q^{-3}(q^4+q^2-1) [n]_2 b_3 b_5^{n-1}-q^{-1}\qmq{1} [n]_2 b_2 b_4 b_5^{n-1},\\
 b_5^n b_2 = b_2 b_5^n-\qmq{1} q^{2-3 n} [n]_2 b_4^2 b_5^{n-1},\qquad
 b_5^n b_4 = q^{-3 n} b_4 b_5^n,\\
 b_4^n b_1 = -\qmq{3} q^{6-3 n} [n-2]_1 [n-1]_1 [n]_1 b_3^2 b_4^{n-3}-\qmq{1} q^{-n} [n]_1b_2^2 b_4^{n-1}\\
\hphantom{b_4^n b_1 =}{}
 -\qmq{3} q^{1-2 n} [n-1]_1 [n]_1 b_2 b_3 b_4^{n-2}+b_1 b_4^n,\\
 b_4^n b_2 = [3]_1 q^{2-2 n} [n]_1 b_3 b_4^{n-1}+q^{-n} b_2 b_4^n,\qquad
 b_4^n b_3 = q^{-3 n} b_3 b_4^n,\\
 b_3^n b_1 = q^{-3 n} b_1 b_3^n+\qmq{1}^2 q^{3-6 n} [n]_2 [3]_1^{-1}b_2^3 b_3^{n-1},\qquad
 b_3^n b_2 = q^{-3 n} b_2 b_3^n,\qquad
 b_2^n b_1 = q^{-3 n} b_1 b_2^n.
\end{gather*}
The lemma is a direct consequence of these formulae.
\end{proof}

A part of the above results have also been obtained in \cite{Xi95}.

Let us turn to the representations $\pi_{{\bf i}}$ of $A_q(G_2)$.
The elements $\sigma_i$ in Def\/inition~\ref{def:sigma}
and $\sigma_ie_i$
are given by
\begin{gather}\label{g:kobetsu}
\sigma_1 = t_{17}, \qquad
\sigma_2 = t_{26}t_{17}-qt_{27}t_{16},\qquad
\sigma_1 e_1 = t_{27},\qquad
\sigma_2e_2 = t_{36}t_{17}-qt_{37}t_{16}.
\end{gather}
From (\ref{a:cpro}) and (\ref{a:piTg}) with $\alpha_i=\mu_i=1$, we have
\begin{gather*}
\pi_{\bf 1}(\sigma_1)=\bk_1 \bK_2 \bk^{2}_3 \bK_4 \bk_5,\\
\pi_{\bf 1}(\sigma_2)=\bK_2 \bk^{3}_3 \bK^{2}_4 \bk^{3}_5 \bK_6,\\
\pi_{\bf 1}(\sigma_1e_1)={\bf a}^+_1 \bK_2 \bk^{2}_3 \bK_4 \bk_5,\\
\pi_{\bf 1}(\sigma_2e_2)=\bk^{3}_1 \bK^{2}_2 \bk^{3}_3 \bK_4 {\bf A}^+_6
 + [2]_2 \bk^{3}_1 {\bf A}^-_2 \bK_2 {\bf A}^+_4 \bK_4 \bk^{3}_5 \bK_6
 + \am_1^3 {\bf A}^+_2 \bk^{3}_3 \bK^{2}_4 \bk^{3}_5 \bK_6\\
 \hphantom{\pi_{\bf 1}(\sigma_2e_2)=}{}
 + [3]_1 {\bf a}_1^{-2} \bk_1 {\bf a}^+_3 \bk^{2}_3 \bK^{2}_4 \bk^{3}_5 \bK_6
 + [3]_1 {\bf a}^-_1 \bk^{2}_1 \bK_2 \bk^{2}_3 \bK_4 {\bf a}^+_5 \bk^{2}_5 \bK_6\\
 \hphantom{\pi_{\bf 1}(\sigma_2e_2)=}{}
 - q [3]_1 \bk^{3}_1 {\bf A}^-_2 \bK_2 \bk^{2}_3 {\bf A}^+_4 \bK_4 \bk^{3}_5 \bK_6
 + [3]_1 \bk^{3}_1 \bK^{2}_2 {\bf a}^-_3 \bk^{2}_3 {\bf a}_5^{+2} \bk_5 \bK_6
 + \bk^{3}_1 \bK^{2}_2 \bk^{3}_3 {\bf A}^-_4 {\bf a}_5^{+3} \bK_6\\
 \hphantom{\pi_{\bf 1}(\sigma_2e_2)=}{}
 + [3]_1 {\bf a}^-_1 \bk^{2}_1 {\bf A}^-_2 {\bf a}_3^{+2} \bk_3 \bK^{2}_4 \bk^{3}_5 \bK_6
 + [3]_1 {\bf a}^-_1 \bk^{2}_1 \bK_2 {\bf a}^-_3 \bk_3 {\bf A}^+_4 \bK_4 \bk^{3}_5 \bK_6\\
\hphantom{\pi_{\bf 1}(\sigma_2e_2)=}{}
 + \bk^{3}_1 {\bf A}_2^{-2} {\bf a}_3^{+3} \bK^{2}_4 \bk^{3}_5 \bK_6
 + [3]_1 \bk^{3}_1 {\bf A}^-_2 \bK_2 {\bf a}^+_3 \bk_3 \bK_4 {\bf a}^+_5 \bk^{2}_5 \bK_6\\
\hphantom{\pi_{\bf 1}(\sigma_2e_2)=}{}
 + \bk^{3}_1 \bK^{2}_2 {\bf a}_3^{-3} {\bf A}_4^{+2} \bk^{3}_5 \bK_6
 + [3]_1 \bk^{3}_1 \bK^{2}_2 {\bf a}_3^{-2} \bk_3 {\bf A}^+_4 {\bf a}^+_5 \bk^{2}_5 \bK_6,\\
\lambda_1^{-1}\pi_{\bf 1}(\xi_1)={\bf a}^+_1 \bk^{-1}_1,\\
\lambda_2^{-1}\pi_{\bf 1}(\xi_2)={\bf a}_1^{-3} {\bf A}^+_2 \bK^{-1}_2
 + [2]_2 \bk^{3}_1 {\bf A}^-_2 \bk^{-3}_3 {\bf A}^+_4 \bK^{-1}_4
 - q [3]_1 \bk^{3}_1 {\bf A}^-_2 \bk^{-1}_3 {\bf A}^+_4 \bK^{-1}_4\\
\hphantom{\lambda_2^{-1}\pi_{\bf 1}(\xi_2)=}{}
 + [3]_1 {\bf a}_1^{-2} \bk_1 \bK^{-1}_2 {\bf a}^+_3 \bk^{-1}_3
 + [3]_1 {\bf a}^-_1 \bk^{2}_1 {\bf a}^-_3 \bk^{-2}_3 {\bf A}^+_4 \bK^{-1}_4
 + [3]_1 {\bf a}^-_1 \bk^{2}_1 \bk^{-1}_3 \bK^{-1}_4 {\bf a}^+_5 \bk^{-1}_5\\
\hphantom{\lambda_2^{-1}\pi_{\bf 1}(\xi_2)=}{}
 + \bk^{3}_1 \bK_2 \bK^{-1}_4 \bk^{-3}_5 {\bf A}^+_6 \bK^{-1}_6
 + [3]_1 {\bf a}^-_1 \bk^{2}_1 {\bf A}^-_2 \bK^{-1}_2 {\bf a}_3^{+2} \bk^{-2}_3
\\
\hphantom{\lambda_2^{-1}\pi_{\bf 1}(\xi_2)=}{}
 + [3]_1 \bk^{3}_1 {\bf A}^-_2 {\bf a}^+_3 \bk^{-2}_3 \bK^{-1}_4 {\bf a}^+_5 \bk^{-1}_5
 + \bk^{3}_1 {\bf A}_2^{-2} \bK^{-1}_2 {\bf a}_3^{+3} \bk^{-3}_3
\\
\hphantom{\lambda_2^{-1}\pi_{\bf 1}(\xi_2)=}{}
 + [3]_1 \bk^{3}_1 \bK_2 {\bf a}^-_3 \bk^{-1}_3 \bK^{-2}_4 {\bf a}_5^{+2} \bk^{-2}_5
 + \bk^{3}_1 \bK_2  {\bf A}^-_4 \bK^{-2}_4 {\bf a}_5^{+3} \bk^{-3}_5\\
\hphantom{\lambda_2^{-1}\pi_{\bf 1}(\xi_2)=}{}
 + \bk^{3}_1 \bK_2 {\bf a}_3^{-3} \bk^{-3}_3 {\bf A}_4^{+2} \bK^{-2}_4
 + [3]_1 \bk^{3}_1 \bK_2 {\bf a}_3^{-2} \bk^{-2}_3 {\bf A}^+_4 \bK^{-2}_4 {\bf a}^+_5
\bk^{-1}_5,\\
\pi_{\bf 2}(\sigma_1)=\bk_2 \bK_3 \bk^{2}_4 \bK_5 \bk_6,\\
\pi_{\bf 2}(\sigma_2)=\bK_1 \bk^{3}_2 \bK^{2}_3 \bk^{3}_4 \bK_5,\\
\pi_{\bf 2}(\sigma_1e_1)=\bK_1 \bk^{2}_2 \bK_3 \bk_4 {\bf a}^+_6
 + {\bf A}^-_1 {\bf a}^+_2 \bK_3 \bk^{2}_4 \bK_5 \bk_6
 + \bK_1 \bk^{2}_2 \bK_3 {\bf a}^-_4 {\bf A}^+_5 \bk_6
 + \bK_1 {\bf a}_2^{-2} {\bf A}^+_3 \bk^{2}_4 \bK_5 \bk_6\\
\hphantom{\pi_{\bf 2}(\sigma_1e_1)=}{}
 + [2]_1 \bK_1 {\bf a}^-_2 \bk_2 {\bf a}^+_4 \bk_4 \bK_5 \bk_6
 + \bK_1 \bk^{2}_2 {\bf A}^-_3 {\bf a}_4^{+2} \bK_5 \bk_6,\\
\pi_{\bf 2}(\sigma_2e_2)={\bf A}^+_1 \bk^{3}_2 \bK^{2}_3 \bk^{3}_4 \bK_5,\\
\lambda_1^{-1}\pi_{\bf 2}(\xi_1)={\bf A}^-_1 {\bf a}^+_2 \bk^{-1}_2\!
 + [2]_1 \bK_1 {\bf a}^-_2 \bK^{-1}_3 {\bf a}^+_4 \bk^{-1}_4
 + \bK_1 {\bf a}_2^{-2} \bk^{-1}_2 {\bf A}^+_3 \bK^{-1}_3
 + \bK_1 \bk_2 {\bf a}^-_4 \bk^{-2}_4 {\bf A}^+_5 \bK^{-1}_5\\
\hphantom{\lambda_1^{-1}\pi_{\bf 2}(\xi_1)=}{}
 + \bK_1 \bk_2 \bk^{-1}_4 \bK^{-1}_5 {\bf a}^+_6 \bk^{-1}_6
 + \bK_1 \bk_2 {\bf A}^-_3 \bK^{-1}_3 {\bf a}_4^{+2} \bk^{-2}_4,\\
\lambda_2^{-1}\pi_{\bf 2}(\xi_2)={\bf A}^+_1 \bK^{-1}_1.
\end{gather*}

We f\/ind that $\pi_{\bf i}(\sigma_i)$ is invertible.
Comparing these formulae with Lemma \ref{lemma:rho-rep-G}
by using (\ref{a:ketta}), the equality (\ref{a:key}) is directly checked.
Thus Proposition \ref{a:pr:rnk2} is established for $G_2$.

In terms of the checked intertwiner $\mathscr{F}$ in Section \ref{a:ss:g2},
Theorem \ref{theorem:gamma=R} implies
\begin{gather*}
E^{a,b,c,d,e,f}_{\bf 2} = \sum_{i,j,k,l,m,n}
\mathscr{F}^{a b c d e f}_{i j k l m n} E^{n,m,l,k,j,i}_{\bf 1}.
\end{gather*}

\section{Discussion}\label{a:sec6}

In view of Proposition \ref{a:pr:rnk2} it is natural to expect that the map def\/ined on generators
of $U_q^+(\geh)$ as $e_i \mapsto \eta_i:=\sigma_ie_i/\sigma_i$ extends to an algebra homomorphism from
$U_q^+(\geh)$ to $A_q(\geh)_\mathcal{S}$, namely, $\eta_i$~satisf\/ies $q$-Serre relations. In fact, it is
true not only for rank 2 cases but also for any $\geh$.

\begin{theorem}\label{th:serre}
In $A_q(\geh)_\mathcal{S}$ the following relation holds for any $i$, $j$ $(i\ne j)$:
\[
\sum_{r=0}^{1-a_{ij}}(-1)^r\eta_i^{(r)}\eta_j \eta_i^{(1-a_{ij}-r)}=0.
\]
\end{theorem}

\begin{proof}
By relabeling of Dynkin indices we can assume $i=1$, $j=2$.
Set $\tau_i=\sigma_ie_i$ for $i=1,2$. Then from Proposition \ref{pr:q-comm} we have
\begin{gather}\label{eq1}
\sigma_i \tau_i=q_i \tau_i \sigma_i,\quad i=1,2, \qquad
\sigma_i \tau_j=\tau_j \sigma_i, \quad i,j=1,2; \quad i\neq j.
\end{gather}
Using \eqref{add-mult} with these relations one verif\/ies
\[
\eta_1^{r}\eta_2 \eta_1^{s}=q_1^{(r+s)(r+s-1)/2}\big(\tau_1^r\tau_2\tau_1^s\big)/\big(\sigma_1^r\sigma_2\sigma_1^s\big).
\]
Here we have set $s=1-a_{12}-r$. Recalling that $\sigma_1$ and $\sigma_2$ commute with each other,
we can reduce the claim to showing
\[
Z:=\sum_{r=0}^{1-a_{12}}(-1)^r\tau_1^{(r)}\tau_2 \tau_1^{(s)}=0.
\]

Note that the right (resp. left) weight of $Z$ is $(1-a_{12})(\varpi_1-\alpha_1)+(\varpi_2-\alpha_2)$
(resp.\ $w_0((1-a_{12})\varpi_1+\varpi_2)$).
The two weights are not related by the longest element $w_0 \in W$.
Hence if we show
$f_iZ=Z f_i=0$ for any $i$, we can conclude $Z=0$ by
the remark after Lemma~\ref{le:ff}.
The properties $f_iZ=0$ for any $i$ and $Z f_i=0$ for $i\neq 1,2$ are trivial.

First we show $Z f_2=0$. We have
\[
\big(\tau_1^r \tau_2 \tau_1^s\big) f_2
=\tau_1^r (\tau_2 f_2) \big(\tau_1 k_2^{-1}\big)^s
=\tau_1^r \sigma_2 (\beta \tau_1)^s=\beta^s \tau_1^{r+s} \sigma_2,
\]
where $\beta=q_2^{-\langle h_2,\varpi_1-\alpha_1 \rangle}=q_2^{a_{21}}=q_1^{a_{12}}$ and
we have used~\eqref{eq1}. Hence,
\[
Z f_2
=\left(\sum_{r+s=1-a_{12}} \frac{(-q_1^{-a_{12}})^s}{[r]_1![s]_1!} \right)
(-\tau_1)^{1-a_{12}} \sigma_2=0.
\]
In the last equality we have used the following formula:
\[
\sum_{i=0}^m(-z)^i{m\brack i}=\prod_{j=1}^m\big(1-zq^{2j-m-1}\big),
\]
where ${m\brack i}=[m]!/([i]![m-i]!)$.

Next, we show $Z f_1=0$.
\begin{gather*}
\big(\tau_1^r \tau_2 \tau_1^s\big) f_1
= \sum_{i=1}^r \tau_1^{r-i} \sigma_1 \big(\tau_1 k_1^{-1}\big)^{i-1}\big(\tau_2 k_1^{-1}\big) \big(\tau_1 k_1^{-1}\big)^s
+\tau_1^r \tau_2 \sum_{i=1}^s \tau_1^{s-i} \sigma_1 \big(\tau_1 k_1^{-1}\big)^{i-1}\\
\hphantom{\big(\tau_1^r \tau_2 \tau_1^s\big) f_1}{}
= \sum_{i=1}^r \delta \gamma^{i-1-s}\tau_1^{r-1} \tau_2 \tau_1^{s}\sigma_1
+\sum_{i=1}^s \gamma^{i-1}\tau_1^{r} \tau_2 \tau_1^{s-1}\sigma_1,
\end{gather*}
where constants $\gamma$, $\delta$ are determined by $\sigma_1 \big(\tau_1 k_1^{-1}\big)=\gamma \tau_1 \sigma_1$,
$\sigma_1 \big(\tau_2 k_1^{-1}\big)=\delta \tau_2 \sigma_1$ and hence we have
$\gamma=q_1 q_1^{-\langle h_1, \varpi_1-\alpha_1 \rangle}=q_1^2$,
$\delta=q_1^{-\langle h_1, \varpi_2-\alpha_2 \rangle}=q_1^{a_{12}}$.
Then, we obtain
\begin{gather*}
Z f_1
 =\sum_{r+s=1-a_{12}} \frac{(-1)^r}{[r]_1![s]_1!}
\left(\sum_{i=1}^r \delta \gamma^{i-1+s}\tau_1^{r-1} \tau_2 \tau_1^{s}\sigma_1
+\sum_{i=1}^s \gamma^{i-1}\tau_1^{r} \tau_2 \tau_1^{s-1}\sigma_1\right)\\
\hphantom{Z f_1}{}
 =\sum_{r+s=1-a_{12}} \frac{(-1)^r}{[r]_1![s]_1!}
\left(\delta \gamma^{s}\frac{1-\gamma^r}{1-\gamma}\tau_1^{r-1} \tau_2 \tau_1^{s}\sigma_1
+\frac{1-\gamma^s}{1-\gamma}\tau_1^{r} \tau_2 \tau_1^{s-1}\sigma_1\right)\\
\hphantom{Z f_1}{}
 =\sum_{r+s=1-a_{12}} \left(-(-1)^{r-1}q_1^s\tau_1^{(r-1)} \tau_2 \tau_1^{(s)}
+(-1)^r q_1^{s-1}\tau_1^{(r)} \tau_2 \tau_1^{(s-1)}\right)\sigma_1=0
\end{gather*}
as desired.
\end{proof}

\begin{remark}
The special case $w=w_0$ of \cite[Theorem 3.7]{Yak} gives
Theorem~\ref{th:serre} here.
Moreover \cite[Theorem~3.7]{Yak} also
shows that
$U_q^+(\geh)$ is isomorphic to an explicit subalgebra of $A_q(\geh)_
\mathcal{S}$.
We would like to thank the referee for pointing this out and for
giving helpful comments.
\end{remark}

It will be interesting to investigate it further in the light of the
quantum cluster algebra which has been recognized as a
fundamental structure in the quantized algebra of functions~\cite{GLS}.
The representations via multiplication on PBW bases also play a
fundamental role in the study of the positive principal series representations
and  modular double \cite{Ip}.

In this paper we have not discussed
the analogue of the tetrahedron and 3D ref\/lection equations
for general $\geh$.
However, from our proof of Theorem \ref{theorem:gamma=R},
we expect that the basic constituents are
$R$ and $K$ only, and their compatibility condition
is reduced to the rank 2 cases (\ref{a:teq}) and (\ref{a:3dref}).

\subsection*{Acknowledgments}
The authors thank Ivan C.H.~Ip, Anatol N.~Kirillov,
Toshiki Nakashima and Masatoshi Noumi for communications.
They also thank one of the referees for drawing attention to the referen\-ces~\cite{Joseph, Yak}.
This work is supported by Grants-in-Aid for
Scientif\/ic Research No.~23340007, No.~24540203,
No.~23654007 and No.~21340036 from JSPS.

\pdfbookmark[1]{References}{ref}
\LastPageEnding

\end{document}